\documentclass[12pt,a4paper]{article}

\usepackage[margin=2cm]{geometry}
\usepackage[hidelinks]{hyperref}
\usepackage{graphicx}
\usepackage{amssymb, amsmath, amsfonts, bm, mathrsfs,amsthm}
\usepackage{siunitx}
\usepackage[utf8]{inputenc}
\usepackage{cleveref}

\newcommand{\E}[0]{\mathcal{E}}
\renewcommand{\O}[0]{\mathcal{O}}
\newcommand{\R}[0]{\mathbb{R}}
\newcommand{\eps}[0]{\varepsilon}
\renewcommand{\Re}[0]{\operatorname{Re}}

\newcommand{\erf}[0]{\operatorname{erf}}
\newcommand{\Ghat}[0]{\widehat{\Gamma}}
\newcommand{\U}[0]{\mathcal{U}}
\newcommand{\V}[0]{\mathcal{V}}

\newcommand{\grad}[0]{\nabla}
\newcommand{\lapof}[0]{\nabla^2}
\newcommand{\divof}[0]{\nabla \cdot}

\newcommand{\abs}[1]{\left|#1\right|}
\newcommand{\br}[1]{\left(#1\right)}

\newcommand{\pd}[2]{\frac{\partial #1}{\partial #2}}

\renewcommand{\O}[0]{\mathcal{O}}
\newcommand{\nhat}[0]{\widehat{n}}
\newcommand{\that}[0]{\widehat{t}}
\newcommand{\Kbar}[0]{\overline{K}}
\newcommand{\Kabs}[0]{\abs{K}}
\newcommand{\Khat}[0]{\widehat{K}}

\newcommand{\Jhat}[0]{\widehat{J}}
\newcommand{\sgrad}[0]{\nabla_\bot}
\newcommand{\grads}[0]{\nabla^\bot}
\newcommand{\sdiv}[0]{\nabla_\bot \cdot}
\newcommand{\divs}[0]{\nabla^\bot \cdot}

\newcommand{\ds}[0]{\partial_\sigma}
\newcommand{\dx}[0]{\partial_\xi}

\newcommand{\notvc}[1]{#1}

\title{Improving accuracy of volume penalised fluid-solid interactions}

\author{Eric W. Hester \thanks{University of Sydney School of Mathematics and Statistics, Sydney, NSW 2006, Australia (\href{mailto:eric.hester@sydney.edu.au}{eric.hester@sydney.edu.au}).
	 }
\and Geoffrey M. Vasil \footnotemark[1]
\and Keaton J. Burns 
\thanks{Massachusetts Institute of Technology Department of Physics, Cambridge, MA 02139, USA}
\thanks{Center for Computational Astrophysics, Flatiron Institute, Simons Foundation, New York, NY 10010, USA}}

\newtheorem{remark}{Remark}[section]
\newtheorem{definition}{Definition}[section]

\begin{document}

\maketitle

\begin{abstract}
We analyse and improve the volume-penalty method, a simple and versatile way to model objects in fluid flows.
The volume-penalty method is a kind of fictitious-domain method that approximates no-slip boundary conditions with rapid linear damping inside the object.
The method can then simulate complex, moving objects in general numerical solvers
without specialised algorithms or boundary-conforming grids.
Volume penalisation pays for this simplicity by introducing an equation-level error, the \emph{model error}, that is related to the damping time $\eta \ll 1$.
While the model error has been proven to vanish as the damping time tends to zero, previous work suggests convergence at a slow rate of $\O(\eta^{1/2})$.
The stiffness of the damping implies conventional volume penalisation only achieves first order numerical accuracy.
We analyse the volume-penalty method using multiple-scales matched-asymptotics with a signed-distance coordinate system valid for arbitrary smooth geometries.
We show the dominant model error stems from a displacement length that is proportional to a Reynolds number $\Re$ dependent boundary layer of size $\O(\eta^{1/2}\Re^{-1/2})$.
The relative size of the displacement length and damping time leads to multiple error regimes.
Our key finding derives a simple smoothing prescription for the damping that eliminates the displacement length and reduces the model error to $\O(\eta)$ in all regimes.
This translates to second order numerical accuracy.
We validate our findings in several comprehensive benchmark problems and finally combine Richardson extrapolation of the model error with our correction to further improve convergence to $\O(\eta^{2})$.
\end{abstract}

\section{Introduction}
Fluid-solid interactions are a continuing challenge for numerical modelling. 
The main difficulty results from imposing fluid boundary conditions at the interface of complicated solid shapes. 
In the most interesting cases, the location, shape, or topology of a solid co-evolves with the dynamics of the fluid. 

The most obvious way to handle complex boundary conditions is with ``boundary-conforming" grids or numerical elements. 
While many sophisticated numerical packages take this approach \cite{DevilleHighorderMethodsIncompressible2002,FischerNek5000WebPage2008,KarniadakisSpectralHpElement2005,CantwellNektarOpensourceSpectral2015,BlackburnSemtexSpectralElement2019}, the method is expensive to simulate, and challenging to implement from a software perspective. 
There is a need for simpler schemes that are compatible with existing codes designed to solve standard fluid-dynamical models in simple domains.

This paper examines the volume-penalty method (VPM), which approximates no-slip boundaries with rapid damping inside the fictitious solid interior.
The VPM thereby avoids explicit tracking of the boundary; dynamically updating source terms is comparatively easy.
Fluid-stress boundary integrals are also replaced by volume integrals of the damping.
This simple method arose from general forcing approaches \cite{BriscoliniDevelopmentMaskMethod1989,GoldsteinModelingNoSlipFlow1993},
and is motivated physically by porous-media flows \cite{AngotAnalysisSingularPerturbations1999,AngotPenalizationMethodTake1999}.
Its versatility has spawned applications to simulations of the Navier-Stokes equations with arbitrary geometries  \cite{KhadraFictitiousDomainApproach2000,KevlahanComputationTurbulentFlow2001,IaccarinoImmersedBoundaryTechnique2003,YenApproximationLaplaceStokes2014}, 
moving boundaries \cite{MittalImmersedBoundaryMethods2005,KolomenskiyAnalysisDiscretizationVolume2015}, 
magnetohydrodynamics \cite{SchneiderPseudospectralMethodVolume2011}, 
flux boundary conditions \cite{KadochVolumePenalizationMethod2012,SakuraiVolumePenalizationInhomogeneous2019}, 
and even insect flight \cite{EngelsNumericalSimulationFluid2015,EngelsFluSINovelParallel2016}. 

Philosophically, boundary conditions exist to model rapid media changes. 
A ``solid'' wall is, in reality, an extreme stiffness change in the elastic properties of a continuum.
Students in physics often first learn about boundary conditions for Maxwell's equations with ``pill-box'' arguments applied to rapidly-changing macroscopic media. 
The VPM simply backs off from these notions and allows some level of resolved rapid change in medium properties. 

There is a long history of this approach, disguised in many forms. 
The Rankine-Hugoniot jump conditions emerge in the limit of vanishing diffusivity in Burgers equation \cite{ColeQuasiLinearParabolicEquation1951,HopfPartialDifferentialEquation1950}. 
Gibbs' analysis of capillary forces \cite{GibbsEquilibriumHeterogeneousSubstances1878} and Stefan's melting phase transitions \cite{StefanUeberTheorieEisbildung1891} supposed material discontinuities.
Further studies found small-but-finite length scales standing in for infinitesimal jumps (e.g., Van der Waals' explanation of surface tension \cite{VanderWaalsThermodynamicTheoryCapillarity1979}, and Cahn and Hilliard's model of phase separation \cite{CahnFreeEnergyNonuniform1958}). 

Many methods exploit smoothness to design simple algorithms for simulating PDEs on complicated domains. 
Phase-field models are a prime example of this smooth approach \cite{HohenbergTheoryDynamicCritical1977,LangerModelsPatternFormation1986,BoettingerPhaseFieldSimulationSolidification2002}. 
Diffuse-interface methods extend the phase-field approach to simulate two-phase liquids \cite{JacqminCalculationTwoPhaseNavier1999,AndersonDiffuseInterfaceMethodsFluid1998}. 
The diffuse-domain method also originates from phase-field models and approximates Neumann, Robin, and Dirichlet boundary conditions \cite{LiSolvingPDEsComplex2009}. 
Fictitious-domain methods are similar to diffuse-domain methods and solve for Lagrange multipliers within the object to enforce the desired boundary conditions \cite{GlowinskiFictitiousDomainMethod1994,GlowinskiFictitiousDomainApproach2001,YuDirectforcingFictitiousDomain2007,YuDLMFDMethod2005}. 
The immersed-boundary method is another classical method for simulating elastic objects in fluids \cite{PeskinNumericalAnalysisBlood1977,PeskinImmersedBoundaryMethod2002} (though rigid objects can cause numerical stiffness \cite{MittalImmersedBoundaryMethods2005}).
Of course, the length scales of realistic medium transitions are often too small to resolve computationally.
In these situations, smooth models are a numerical appliance which trade the accuracy of boundary conditions for greater computational simplicity.

The key challenge of the VPM, and other smoothed models, is to minimise the error compared with formal boundary conditions.
Much progress has been made in improving the error of phase-field models \cite{CaginalpAnalysisPhaseField1986,CaginalpStefanHeleShawType1989,KarmaPhasefieldMethodComputationally1996,ChenRapidlyConvergingPhase2006} and diffuse-domain methods \cite{FranzNoteConvergenceAnalysis2012,LervagAnalysisDiffusedomainMethod2015,BurgerAnalysisDiffuseDomain2017,PoulsenSmoothedBoundaryMethod2017,YuHigherorderAccurateDiffusedomain2020}.
Progress has been slower for the VPM.

Previous work has shown that the error of the VPM vanishes in the limit that the damping timescale $\eta$ tends to zero \cite{AngotAnalysisSingularPerturbations1999}.
Later work showed this error converges at the slow rate $\O(\eta^{1/2})$ \cite{CarbouBoundaryLayerPenalization2003}.
The stiffness of the damping implies this error leads to merely first-order numerical convergence.
These works also used discontinuous damping, which restricts numerical convergence, causing Gibbs phenomena in spectral codes \cite{KevlahanComputationTurbulentFlow2001}.

In this paper, we improve upon previous understanding of the VPM and derive simple corrections that reduce the overall error. 
Our analysis extends the higher-order analysis of phase-field models \cite{KarmaPhasefieldMethodComputationally1996,ChenRapidlyConvergingPhase2006} and diffuse-domain methods \cite{PoulsenSmoothedBoundaryMethod2017,YuDirectforcingFictitiousDomain2007}. 
We begin in \cref{sec:definitions} by describing the VPM, and the various metrics used to measure the error.
The central penalty quantity is the non-dimensional damping length scale, $\eps$, which depends on the Reynolds number Re and non-dimensional damping time scale, $\eta$. 
\Cref{sec:asymptotics} then develops a general and straightforward multiple-scales matched-asymptotics procedure in powers of $\eps$. 
We use this analysis to quantify the error for arbitrary smooth objects in three dimensions. 
In \cref{sec:sdf}, we introduce the signed-distance coordinate system. 
This analytical tool helps streamline complicated calculations and unambiguously accounts for perplexities such as interface curvature. 
We analyse the VPM up to second-order and derive simple conditions to eliminate the leading-order error. 
\Cref{sec:inner} develops efficient methods to calculate these optimisations in general cases. 
The remainder of the paper verifies the improvement to performance in several test problems. 
In \cref{sec:unsteady-cylinder}, we test our prescriptions in time-dependent accelerating flow past a rotating 2D cylinder. 
We outline the efficient spectral code Dedalus, use it to generate numerical reference and penalised solutions, and validate the improved accuracy of our corrections. 
We finally combine our optimisations with Richardson extrapolation of the model error to further reduce model error to $\O(\eta^2)$.

\section{The volume-penalty method}
\label{sec:definitions}

\subsection{Objects in incompressible hydrodynamics}
\begin{figure}[ht]
	\centering
    \includegraphics[width=\linewidth]{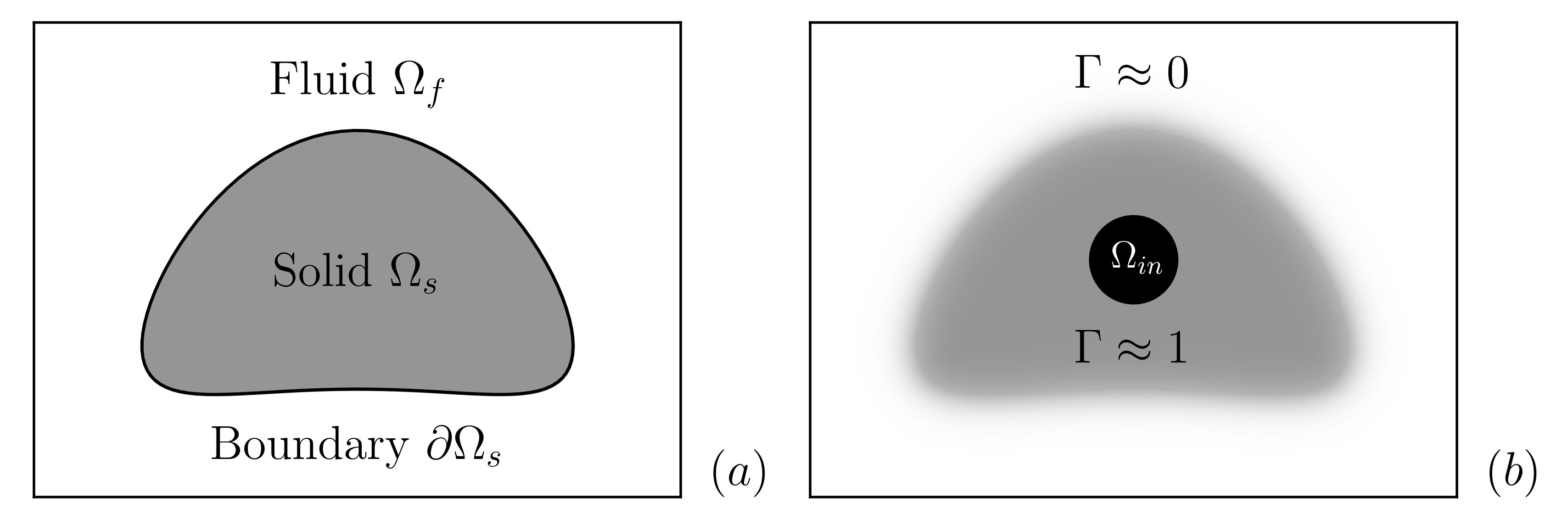}
    \caption{Figure $(a)$ illustrates the standard no-slip boundary formulation of fluid-solid interactions, where the domain is separated into the fluid domain $\Omega_f$, and the solid domain $\Omega_s$.
  	Figure $(b)$ shows the volume penalised approach, in which the solid is represented by a mask function $\Gamma$ which exists throughout the domain. The function need not be discontinuous, or coincide exactly with the intended solid domain. We account for possible internal no-slip boundaries $\Omega_{in}$ in the solid object, which are required in \cref{sec:unsteady-cylinder}.}
    \label{fig:domain-diagrams}
\end{figure}

Fluid-solid interactions in incompressible viscous fluids are normally modelled by partitioning the domain $\Omega$ into solid $\Omega_s$ and fluid $\Omega_f$ components (see \cref{fig:domain-diagrams} $(a)$).
On $\Omega_f$ we solve the Navier-Stokes equations,
	\begin{align}\label{eq:ns-true}
	\partial_t u + \frac{1}{\rho}\nabla{p} - \nu \nabla^2{{u}} + {u}\cdot \nabla{{u}} - {f} &= 0, \qquad \text{and} \qquad
	\nabla \cdot {u} = 0,
	\end{align}
for the fluid velocity ${u}$ and the pressure $p$, 
given a kinematic viscosity $\nu$, constant density $\rho$, and external body force ${f}$.
The solid then acts through no-slip boundary conditions 
	\begin{align}\label{eq:no-slip}
	{u} = {u}_s \quad\text{on}\quad\partial \Omega_s,
\end{align}
where ${u}_s$ is the solid velocity at the boundary.
The fluid in turn influences the solid through drag ${F}$ and torque ${T}$ surface stress integrals
    \begin{align}\label{eq:surface-stress}
    {F} &= \int_{\partial \Omega_s} \Sigma \cdot {\nhat}\, dA, \qquad \text{and} \qquad
    {T} = \int_{\partial \Omega_s} {r}\times (\Sigma \cdot {\nhat})\, dA,
    \end{align}
where $\nhat$ is the outward unit normal vector on the solid boundary $\partial \Omega_s$, and $\Sigma$ is the stress
    \begin{align*}
    \Sigma = -p I + \rho\nu (\nabla {u} + \nabla {u}^{\top}).
    \end{align*}
Unfortunately this formulation can be difficult to discretise in elaborate geometries.
This motivates the volume-penalty method, an alternate model of fluid-solid interactions which allows simpler discretisation, at the cost of some error in the mathematical solution.

\subsection{The volume-penalty approximation}
The volume-penalty method models objects in fluids by replacing no-slip boundary conditions with an added linear damping in the momentum equation,
	\begin{align}\label{eq:ns-penalised}
	\partial_t {u} + \frac{1}{\rho}\nabla{p} - \nu \nabla^2{{u}} + {u} \cdot \nabla{{u} } - {f} &= -\frac{1}{\tau} \Gamma(x,t) ({u} - {u}_s),
	\end{align}
where the \emph{mask function} $\Gamma : \Omega \to [0,1]$ represents the location of the solid (which may evolve in time).
The penalty term has a physical interpretation.
It corresponds to the Darcy drag of the Brinkman equations for fluid flow in a porous medium \cite{AngotPenalizationMethodTake1999}.
The small damping timescale $\tau$ forces the fluid velocity to tend to the solid velocity within the object.
Traditionally, the penalty mask $\Gamma$ is chosen to be a discontinuous indicator function over the solid region \cite{AngotAnalysisSingularPerturbations1999,AngotPenalizationMethodTake1999,KolomenskiyFourierSpectralMethod2009,SchneiderPseudospectralMethodVolume2011,YenApproximationLaplaceStokes2014,SchneiderImmersedBoundaryMethods2015,EngelsNumericalSimulationFluid2015,EngelsFluSINovelParallel2016}.
However, we emphasise that this is not the only choice.
We will later show that it is in fact suboptimal.

The volume-penalty method also approximates the force $F$ and torque $T$ surface integrals (\cref{eq:surface-stress}, \cite{BriscoliniDevelopmentMaskMethod1989}) with simpler-to-discretise volume-penalty integrals \cite{AngotAnalysisSingularPerturbations1999}
	\begin{align}\label{eq:standard-force-calculation}
	{F} &= \rho\int_\Omega \frac{\Gamma}{\tau} \,({u} - {u}_s)\, dV +\rho \int_{\Omega_s} \,(\partial_t u_s - {f}) \, dV + \int_{\partial\Omega_{in}} \Sigma \cdot \nhat \, dS,\\
	{T} &= \rho\int_\Omega \frac{\Gamma}{\tau} \, r \times ({u} - {u}_s)\, dV + \rho\int_{\Omega_s} {r} \times (\partial_t u_s - {f}) \, dV + \int_{\partial\Omega_{in}} \! r \times (\Sigma \cdot \nhat) \, dS.
	\end{align}
The second terms account for forcing within the object $\Omega_s$ and are analogous to the unsteady correction from \cite{EngelsNumericalSimulationFluid2015,UhlmannImmersedBoundaryMethod2005}.
The final corrections account for possible no-slip boundaries $\partial\Omega_{in}$ {within} the volume penalised object (as in \cref{fig:domain-diagrams}).
Both corrections are necessary in \cref{sec:unsteady-cylinder}.

The volume-penalty method is useful because it removes the object boundary from the mathematical formulation.
A single set of equations applies throughout the entire domain,
the solid is represented by a mask function $\Gamma$,
and the force and torque are calculated with volume integrals.
This makes it simple to implement numerically,
as we do not need the numerical grid to conform to the solid, nor do we need to interpolate quantities at the boundary.

However, the volume-penalty model is different to the standard no-slip model.
The mathematical solution of each problem will therefore differ.
We define this equation-level difference as the \emph{model error}.

\begin{definition}[Model error]\label{def:model-error}
The difference between the mathematical solution of the true no-slip boundary formulation, and the volume-penalty approximation. This excludes numerical differences due to discretisation.
\end{definition}

We emphasise that the model error exists separate to the \emph{discretisation  error},
which is specific to a particular numerical implementation.
The model error is {implementation-independent}, and therefore more fundamental.
By finding mathematical ways to minimise the model error, we derive improvements that apply to all possible implementations of the volume-penalty method.

But what \emph{is} the model error?
And what should it depend on?
We answer these questions in the remainder of this section.

\subsection{Non-dimensional penalty parameters}
The most important determinants of the model error are two non-dimensional penalty parameters. 
We initiate our non-dimensionalisation by defining the dimensional damping time $\tau$, the characteristic length scale of the domain $L$, and the characteristic velocity scale of the fluid $U$.
These induce three key non-dimensional parameters, the \emph{Reynolds number} $\Re$, the \emph{damping time scale} $\eta$, and the \emph{damping length scale} $\eps$
    \begin{align}
    \Re &\equiv \frac{U L}{\nu}, & \eta &\equiv \frac{\tau\,U}{L}, & \eps &\equiv \frac{\sqrt{\nu \tau}}{L}.
    \end{align}
The Reynolds number $\Re$ characterises the ratio of viscous and inertial time scales,
the (non-dimensional) damping \emph{time} scale $\eta$ measures the separation of damping and inertial time scales, 
and the damping \emph{length} scale $\eps$ comes from balancing the viscous and damping terms.

The damping time $\eta$ and length $\eps$ are closely related via the Reynolds number $\Re$
    \begin{align}
    {\eta} = {\Re  \eps^2},
    \end{align}
and are {equal} when $\eta = \Re^{-1} = \eps$.
The length and time scale each induce different kinds of model errors.
Either side of this equality one scale dominates, leading to distinct error regimes.

We then write the penalised equations as
\begin{align}
\partial_t{u} + u \cdot \nabla{u} + \nabla{p} - \frac{1}{\Re}\nabla^2{u} - f &= - \frac{1}{\eta} \Gamma(x,t) (u - u_s),
\end{align}
where $u$ and $p$ are penalised solutions.
Solutions to the `true' non-dimensional Navier-Stokes equations with no-slip boundary conditions are denoted with a zero subscript (e.g. $u_0$, $p_0$), 
in keeping with the convergence of the penalised solution $u$ to the true solution $u_0$ as $\eps \to 0$.

\begin{remark}[Error vs {computational effort}]
To numerically resolve the damping length $\eps$ and time scale $\eta = \Re\eps^2$,
the spatial degrees of freedom (Fourier modes, grid points etc.) must scale as $\eps^{-D}$ in dimension $D$,
and total time steps must scale as $\eps^{-2}$,
for total effort $\eps^{-(D + 2)}$.
Given an error scaling of $\eps^{\alpha}$, the increased effort to halve the error is proportional to $2^{(D + 2)/\alpha}$.
This implies a 16-fold increase in cost in 2D and 32-fold increase in 3D for the existing scaling $\alpha=1$.
It is thus {essential} to improve the error scaling $\alpha$.	
\end{remark}

\subsection{Mask functions}
The final quantity that influences the model error is the choice of mask function $\Gamma$.
We now define normalised mask functions $\Ghat$, and the general masks $\Gamma$ that we consider in this paper.
\begin{definition}[Normalised mask function $\Ghat$]
A function on $\R$ which satisfies
	\begin{enumerate}
	\item Boundedness:	$\Ghat: \R \to [0, 1]$.
	\item Limiting behaviour: $\lim_{x\to-\infty}\Ghat(x) = 1$ and $\lim_{x\to\infty}\Ghat(x) = 0$.
	\item Monotonicity: $x_2 > x_1 \implies \Ghat(x_2) < \Ghat(x_1)$.
    \item Symmetry: $\Ghat(x) + \Ghat(-x) = 1$.
    \item Normalisation: $\frac{d\Ghat}{dx}\big|_{x=0} = -1$.
    \end{enumerate}
\end{definition}
We concentrate on smooth normalised mask functions and more restrictive compact normalised mask functions.
Any smooth normalised mask function $\Ghat$ can be compactified to the interval $[-c,c]$ using the following transformation \cite{BoydFourierEmbeddedDomain2005},
	\begin{align}
	\Ghat^c(x) = \Ghat \Big( \tfrac{x}{\sqrt{1-{x^2}/{c^2}}}\Big) \quad \text{where} \quad -c < x < c.
	\end{align}
General mask functions $\Gamma_{\ell,\delta}$ can be generated from a normalised mask function $\Ghat$ through scaling $\delta$ or shifting $\ell$,
	\begin{align}\label{eq:mask-scaled}
		\Gamma_{\ell,\delta}(x) &= \Ghat\left(\tfrac{x-\ell}{\delta}\right).
	\end{align}	
Discontinuous masks are the limiting case of smooth functions as the thickness approaches zero.	
The standard approach takes the limit $\delta \to 0, \ell = 0$ {first}, then considers convergence as $\eps \to 0$.
This is not optimal.
Choosing a particular $\ell$ or $\delta \propto \eps$ cancels the leading order error of the volume-penalty method.

For concreteness, we define the \emph{hyperbolic tangent mask} $\Ghat_{\tanh}$ and \emph{error function mask} $\Ghat_{\erf}$,
		\begin{align}\label{eq:mask-tanh-erf}
	    \Ghat_{\tanh}(x) &= \frac{1}{2} \left(1 - \tanh 2 x\right), &
		\Ghat_{\erf}(x) &= \frac{1}{2} \left(1 - \erf{\sqrt{\pi} x}\right),
		\end{align}
and compactified versions thereof, $\Khat_{[\tanh;c]}, \Khat_{[\erf;c]}$, on the interval $[-c,c]$.

These definitions can be extended to higher dimensions by defining a locally conformal coordinate system at the mask boundary and applying the one dimensional mask functions in the wall-normal direction.
Before developing this coordinate system in \cref{sec:sdf}, we now give several metrics of the model error.

\subsection{Model error metrics and the displacement length}
There are many possible ways to quantify the model error of the volume-penalty method.
We give two metrics that measure the accuracy in an average and worst-case sense, 
as well as the physically meaningful error in the drag calculation.

\begin{definition}[Global error $\E_1$]\label{def:local-error}
The $L^1$ norm of the difference of the penalised solution $u$ from the reference solution $u_0$ in the fluid domain,
			\begin{equation}\label{eq:local-error}
			\E_1 = ||u-u_0||_{L^1(\Omega_f)} = \frac{1}{|\Omega_f|}\int_{\Omega_f} |u - u_0|\,dV.
			\end{equation}
\end{definition}
\begin{definition}[Local error $\E_\infty$]\label{def:global-error}
		The $L^\infty$ norm of the difference of the penalised solution $u$ from the reference solution $u_0$ in the fluid domain,
			\begin{equation}\label{eq:global-error}
			\E_\infty = ||u-u_0||_{L^\infty(\Omega_f)} = \max_{x\in\Omega_f} |u - u_0|.
			\end{equation}	
\end{definition} 

\begin{definition}[{Drag error} $\Delta F$, Torque error $\Delta T$]\label{def:drag-error}
		The difference between the volume integral for the penalised drag $F$ (torque $T$) and the surface stress integral for the reference drag $F_0$ (reference torque $T_0$) on a true no slip boundary.
	\begin{align}
    \Delta F &= F - F_0,  &     \Delta T &= T - T_0.
    \end{align}
\end{definition}
We also justify the fluid error norms by noting that together they quantify error localisation.
If $\abs{u - u_0} \approx \eps$ in a region of size $d$, and much smaller elsewhere,
then $\E_1 \approx \eps d$ (for a unit volume domain), while $\E_\infty \approx \eps$, 
so the ratio $\E_1/\E_\infty \approx d$ quantifies the degree of error localisation.
		
The key insight of this paper is that conventional masks
do {not} optimise the global error $\E_1$ or the drag error $\Delta F$ with respect to the penalty parameter $\eps$.
This is because in general the far field behaviour of the penalised velocity $u$ is shifted by an amount proportional to the damping length scale $\eps$.
We call this difference the \emph{displacement length} $\ell^*$,
\begin{definition}[{Displacement length} $\ell^*$]
	(i). The difference in size between the desired no-slip boundary $\partial \Omega_s$,
	and a hypothetical no-slip boundary $\partial \Omega_s^*$ the mask function $\Gamma$ ``most closely'' approximates. 
	(ii). Alternatively, the difference in size between the desired no-slip boundary $\partial \Omega_s$,
	and the optimal mask function $\Gamma^*$ that most closely approximates it.
\end{definition}
This definition of the displacement length can refer to the size of the mask or the desired no-slip boundary.
We use these meanings interchangeably as they relate to the same phenomenon:
the ideal mask is in general {not} the same shape or size as the solid we wish to approximate.
The concrete advance of this paper is a simple correction to offset the displacement length error by refining the edge of the mask.
This {localises} the error to the boundary layer, reducing the average fluid error $\E_1$ and the drag and torque errors $\Delta F$ and $\Delta T$ to $\O(\eps^2)$.
Having defined the appropriate concepts by which to measure the model error, we now analyse it using a general multiple-scales matched-asymptotics framework.

\section{Multiple-scales matched-asymptotics analysis of the model error}
\label{sec:asymptotics}
\begin{figure}
	\centering
    \includegraphics[width=\linewidth]{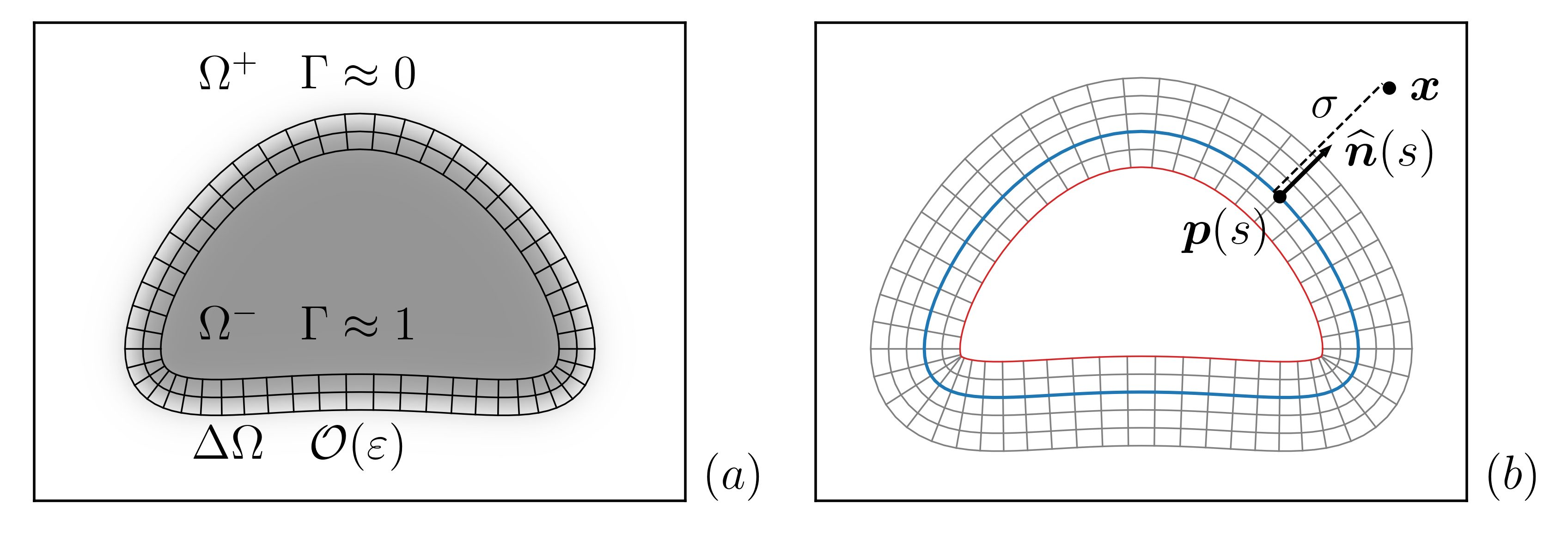}
    \caption{Figure $(a)$ shows the asymptotic regions considered in \cref{sec:asymptotics} for more general masks. The asymptotic fluid domain $\Omega^+$ and solid domain $\Omega^-$ are separated by an $\O(\eps)$ size interfacial region $\Delta \Omega$.
    Figure $(b)$ illustrates the signed distance function coordinate system used in $\Delta \Omega$. 
    Points on the interface $\notvc{p}(s)$ (blue) are parameterised by surface coordinates $s$.
    Points away from the interface $x$ can be reached by moving a distance $\sigma$ in the normal direction $\notvc{\nhat}(s)$ from the closest point on the interface $\notvc{p}(s)$. 
    We also plot level sets of the signed distance $\sigma$ and surface coordinates $s$.
    Coordinate singularities (the corners of the red curve) will occur, but the coordinate system remains well-behaved in the interface region for smooth objects.
    The figure is in two dimensions for clarity, though the analysis of \cref{sec:asymptotics} is done in three dimensions.}
    \label{fig:sdf}
\end{figure}

Understanding volume penalisation requires analysis of the model error.
This analysis must be performed at the equation level,
and should be a straightforward procedure that works in general geometries and flow configurations.
We provide such a framework using multiple-scales matched-asymptotics valid for general tensor-valued equations in arbitrary smooth geometries.
This procedure breaks into several steps:
\begin{enumerate}
	\item Partition the domain into solid $\Omega^-$, fluid $\Omega^+$, and interfacial regions $\Delta \Omega$ (\cref{fig:sdf} $(a)$).
	\item Use the signed distance $\sigma$ as a coordinate in the interfacial region (\cref{sec:sdf}, \cref{fig:sdf} $(b)$).
	\item Rescale the normal coordinate and operators by $\eps$ in the interfacial region (\cref{sec:scaled-sdf}).
	\item Expand the variables in a formal asymptotic power series in $\eps$ each region (\cref{sec:power-series}).
	\item Connect regions with and asymptotic matching conditions (\cref{sec:matching}).
	\item Iterate to solve the zeroth, first, and second order problems (\cref{sec:zeroth,sec:first,sec:second}).
\end{enumerate}
The philosophy of the approach is to determine the asymptotic evolution of the volume-penalty method accurate to second order in $\eps$.
We show that it is possible to analytically determine optimal damping and mask parameters which eliminate the order $\eps$ model error of the volume-penalty method.
We perform the analysis in a body-centred frame incorporating possible acceleration terms into $f$.

The first step partitions the domain into the fluid $\Omega^+$, solid $\Omega^-$, and interfacial regions $\Delta \Omega$, distinguished by the behaviour of the mask function $\Gamma$.
Understanding the interfacial region calls for an appropriate choice of coordinates.

\subsection{Signed-distance coordinates}
\label{sec:sdf}

To understand the volume-penalty method we must analyse the boundary region.
An ideal coordinate system for boundary layers uses the signed distance $\sigma$ of each point to the closest point on the surface.
This choice conforms to the object, aligning the basis vectors with the boundary,
and is {orthogonal}, simplifying the vector calculus.

The signed distance $\sigma$ is defined as the distance of a point $x$ to its nearest point on the interface $p$.
It is not difficult to show the point $x$ must also lie in the direction of the unit normal vector $\nhat$ from point $p$, which we label with surface coordinates $s$ (\cref{fig:sdf} $(b)$),
    \begin{align}\label{eq:sdf}
    x = p(s) + \sigma\, \nhat(s).
    \end{align}
The chain rule allows us to rewrite the normal and surface partial derivatives in terms of the gradient operator, $\partial/\partial \sigma = \nhat \cdot \nabla$, $\partial / \partial s_i = (\partial p/
{\partial}s_i ) \cdot (I + \sigma \left.\nabla\nhat\right|_{\sigma=0}) \cdot \nabla$.
While the normal derivative is straightforward, the tangential derivatives require an understanding of surface differential geometry.
Only then can we invert to find the gradient in the new coordinates.

The first geometric quantity induced by our choice of surface coordinates is a tangent vector basis $t_i$.
Presuming orthogonal surface coordinates, we also derive the dual vector basis $\nabla s_i$ and a unique orthonormal tangent basis $\that_i$,
    \begin{align*}
    {t}_i &= \frac{\partial {p}}{\partial s_i}, &
    \that_i &= \frac{t_i}{|t_i|}, &
    \nabla s_i &= \frac{\that_i}{|t_i|}.
    \end{align*}
This gives us the surface gradient $\sgrad$.
Integrating by parts using the surface area measure $dA = |t_1||t_2|$ gives us the surface divergence $\sgrad \cdot $ of a tangent vector field $u_\bot = u_1 \that_1 + u_2 \that_2$,
    \begin{align*}
    \sgrad &= \nabla s_1  \frac{\partial}{\partial s_1} + \nabla s_2  \frac{\partial}{\partial s_2} = \that_1 \nabla_1 + \that_2 \nabla_2, &
    \sdiv u_\bot &= \frac{1}{|t_1||t_2|}\br{\pd{}{s_1}(|t_2| u_1) + \pd{}{s_2}(|t_1| u_2)}.
    \end{align*}
It is not difficult to show the normal $\nhat$ is everywhere equal to the signed distance gradient.
This implies the gradient of the normal is symmetric, and hence diagonalisable.
The eigenvectors are the principal directions of curvature (and must align with orthogonal surface coordinates), 
and the eigenvalues are the principal curvatures $\kappa_i$,
    \begin{align*}
    \nhat &= \nabla \sigma, &
    \nabla \nhat &= - \kappa_1 \that_1 \that_1 - \kappa_2 \that_2 \that_2 = - K.
    \end{align*}
With an orthonormal frame defined throughout the interfacial region, we decompose the fluid velocity into normal and tangential directions,
    \begin{align*}
   \notvc{u}= u_\sigma \nhat + u_\bot, \quad \text{where} \quad u_\bot = u_1 \that_1 + u_2 \that_2.
    \end{align*}
We can now invert all partial derivative operators and re-express the gradient
    \begin{align}
    \nabla &= \nhat \ds + J^{-1} \cdot \nabla_\bot, \quad \text{where} \quad     J = I - \sigma K,
    \end{align}
where we have defined the scale factor $J$, which possesses several useful identities relating to the mean curvature $\Kbar$, the Gaussian curvature $\Kabs$, and volume measure $dV$, 
    \begin{align*}
    \Kbar &= \kappa_1 + \kappa_2, &
    \Kabs &= \kappa_1 \kappa_2, &
    |J| &= 1 - \sigma \Kbar + \sigma^2 \Kabs, &
    dV = |J| \, d\sigma \, dA.
    \end{align*}
Using the adjugate tensors $\Jhat$ and $\Khat$, we determine the divergence and scalar Laplacian,
    \begin{align*}
    \Khat &= |K| K^{-1} = \kappa_2 \that_1 \that_1 + \kappa_1 \that_2 \that_2,&
	\Jhat &= |J|J^{-1} = I - \sigma \Khat, \\
    \divof{u} &= \frac{\ds (|J| u_\sigma)}{|J|} + \frac{\nabla_\bot \cdot (\Jhat u_\bot)}{|J|}, &
    \lapof f &= \frac{\ds(|J|\ds f)}{|J|} + \frac{\sdiv(\Jhat J^{-1} \sgrad f)}{|J|}.
    \end{align*}
By defining the cross product with the unit normal, which admits simple identities,
    \begin{align*}
    \nabla^\bot &= \nhat \times \nabla_\bot, & 
	    \nabla^\bot \cdot \nabla_\bot &= \nabla_\bot \cdot \nabla^\bot = 0, & 
	    \nhat \times (J u_\bot) &= \Jhat u^\bot, \\
	u^\bot &= \nhat \times u_\bot, &
    	u^\bot \cdot u_\bot &= u_\bot \cdot u^\bot = 0, &
    	\nhat \times u^\bot &= -u_\bot,
    \end{align*}
we can take the curl, and hence find the divergence-free vector Laplacian,
    \begin{align*}
    \nabla \times\notvc{u}&= -\nhat \frac{\sdiv(\Jhat u^\bot)}{|J|} + \Jhat^{-1} (\ds (\Jhat u^\bot) - \grads u_\sigma),\\
    -\nabla \times \nabla\notvc{u}&= 
    \frac{\nhat}{|J|} \br{ - \sdiv(\Jhat J^{-1}(\ds(J u_\bot) - \sgrad u_\sigma )}\\
    &+ \Jhat^{-1} \ds\br{\Jhat J^{-1} (\ds(J u_\bot) - \sgrad u_\sigma)} + \Jhat^{-1} \grads\br{\frac{\divs(J u_\bot)}{|J|}}.
    \end{align*}
We next find the gradient of the basis vectors, and define Ricci rotation coefficients,
    \begin{align*}
    \sgrad \nhat &= - \kappa_1 \that_1\that_1 -\kappa_2 \that_2 \that_2, & 
    &&
    \sgrad \that_i &= \kappa_i \that_i \nhat + \mathcal{R}^{j k}_i \that_j \that_k,\\
    \mathcal{R}^{12}_1 &= \that_1 \cdot (\nabla \that_1) \cdot \that_2 =  \omega_1, &
    \mathcal{R}^{jk}_i &= - \mathcal{R}^{ji}_k, &
    \mathcal{R}^{21}_2 &= \that_2 \cdot (\nabla \that_2) \cdot \that_1 = -\omega_2, &
    \end{align*}
which finally gives us the gradient and advective derivatives of general vector fields   
    \begin{align*}
    \grad\notvc{u}&= \nhat \,\nhat \,\ds u_\sigma + \nhat \ds u_\bot + J^{-1}(\sgrad u_\sigma+ K u_\bot)\,\nhat + J^{-1} (\sgrad u_\bot - K u_\sigma),\\
	u\cdot \nabla f &= u_\sigma \ds f + u_\bot \cdot J^{-1} \sgrad f,\\
	u \cdot \grad\notvc{u}&= \nhat \,\big(u_\sigma \ds u_\sigma+ u_\bot \cdot J^{-1} (\sgrad u_\sigma+ K u_\bot)\big) + u_\sigma \ds u_\bot + u_\bot \cdot J^{-1} (\sgrad u_\bot - K u_\sigma).
    \end{align*}
    
\subsection{Rescaling interfacial coordinates}
\label{sec:scaled-sdf}
The interfacial region is defined as the size $\eps$ region in which the fluid transitions to the solid.
Analysing this region requires rescaling the normal coordinate by $\eps$.
This in turn projects each differential operator into a formal asymptotic series as follows,
    \begin{align*}
    \sigma &= \eps \xi, &
	    J &= 1 - \eps\xi K, & 
	    |J| &= 1 - \eps \xi \Kbar + \eps^2 \xi^2 |K|,\\
    \ds &= \frac{1}{\eps}\dx, &
    J^{-1} &= \sum\nolimits_{k=0}^\infty \eps^k \xi^k K^k, &  
    |J|^{-1} &= \sum\nolimits_{k=0}^\infty \eps^k \xi^k \br{\sum\nolimits_{l=0}^k \kappa_1^l \kappa_2^{k-l}},
   	\end{align*}
	\begin{align*}
    \nabla 	&= \eps^{-1} \nhat \dx + \sum\nolimits_{k=0}^\infty \eps^k \xi^k K^k \sgrad,\\
    |J|\divof\notvc{u}&= \eps^{-1} {\dx u_\sigma}  -\dx(\xi \Kbar u_\sigma) + \sdiv{u_\bot} + \eps \br{\dx(\xi^2 \Kabs u_\sigma) - \xi \sdiv(\Khat u_\bot)},\\
   \notvc{u}\cdot \grad f &= \eps^{-1} u_\sigma \dx f + \sum\nolimits_{k=0}^\infty \eps^k \xi^k K^k u_\bot\cdot \sgrad f\\
	\notvc{u} \cdot \grad\notvc{u}&=\eps^{-1}\left( u_\sigma \dx u_\bot + \nhat u_\sigma\dx u_\sigma\right)\\
	&\quad + \sum\nolimits_{k=0}^{\infty} \eps^k\xi^k u_\bot K^k \cdot \br{(\sgrad u_\bot - K u_\sigma) + (\sgrad u_\sigma+ K u_\bot)\nhat} \nonumber,\\
	\lapof f &= \eps^{-2} \dx^2 f + \eps^{-1}(-\Kbar \dx f) + \eps^0(-\xi \overline{K^2}\dx f + \sgrad \cdot \sgrad f) + \O(\eps)\\
    -\nabla\times\nabla\times\notvc{u}&= \eps^{-2} \dx^2 u_\bot + \eps^{-1} \br{- \Kbar \dx u_\bot - \dx \sgrad u_\sigma- \nhat \sgrad \cdot (\dx u_\bot)} + \O(\eps^0)
    \end{align*}

\subsection{Variable expansions with formal power series}
\label{sec:power-series}

After splitting the domain into the fluid $\Omega^+$, solid $\Omega^-$, and interfacial $\Delta \Omega$ regions, 
each variable $q$ in each region (fluid $q^+$, solid $q^-$, and interfacial $q$) is expressed as a power series in $\eps$,
	\begin{align}
	q^+(x,t) &=    \sum_{k=0}^\infty \eps^k q^+_k(x,t), &
	q^-(x,t) &=    \sum_{k=0}^\infty \eps^k q^-_k(x,t), &
	q(\xi,s,\tau) &= \sum_{k=0}^\infty \eps^k q_k(\xi,s,\tau).
	\end{align}
We perform this expansion for the fluid velocity $u = u_\sigma \nhat + u_\bot$ and pressure $p$.
We then substitute these variable expansions into the infinite hierarchy of equations generated by \cref{sec:scaled-sdf} to derive a system of equations at each order of $\eps$.
Collecting terms and solving the model equations at each order requires a matching procedure between adjacent regions, which we describe in \cref{sec:matching}.

\subsection{Asymptotic matching}
\label{sec:matching}
To ensure agreement between different regions we specify asymptotic matching boundary conditions.
This subtle notion requires asymptotic agreement in intermediate zones $\xi \sim \eps^{-1/2}$ in the limit that $\eps\to 0$.
We can then let $\xi$ approach infinity for the inner problem without encountering coordinate singularities (provided $ \eps \ll \min_i|\kappa_i^{-1}|$), 
and let $\sigma$ approach zero for the outer problem without entering the interfacial region.
That is, for any variable $q$, we require
    \begin{align*}
	\lim_{\eps\to 0} q(\pm\eps^{-1/2}\xi,s,t) &\sim \lim_{\eps\to 0} q^{\pm}(p(s,t)\pm\eps^{+1/2}\,\xi \,\nhat(s,t),t).
    \end{align*}
Each variable is expressed as an asymptotic series.
Each series term in the outer variables can be further expanded as a Taylor series about the interface, involving normal derivatives $\partial_\sigma$ and the rescaled signed distance $\eps \xi = \sigma$.
Combining these expansions and solving order by order in $\eps$, the matching conditions at each order simplify to
    \begin{align*}
    \lim_{\xi\to\pm\infty} q(\xi) = \lim_{\xi\to\pm\infty} \sum_{k=0} \eps^k q_k 
    \sim \lim_{\xi\to\pm\infty}  \sum_{k=0}\eps^k\br{\sum_{\ell=0}^k\frac{\xi^{\ell}}{\ell!} \ds^\ell q^{\pm}_{k-\ell}(0^{\pm})}.
    \end{align*}
Written explicitly, we expect the limiting behaviour of the inner expansions to be asymptotically equivalent to interfacial boundary conditions on the outer variables as $\xi \to \pm\infty$,
\newcommand{\at}[1]{\left.#1\right|}
	\begin{align*}
	\lim_{\xi\to\pm\infty}q_0(\xi,s,t)&\sim \at{q^\pm_0}_{\sigma=0},\\
	\lim_{\xi\to\pm\infty}q_1(\xi,s,t)&\sim \at{\ds q^\pm_0}_{\sigma=0} \xi + \at{q^\pm_1}_{\sigma=0},\\
	\lim_{\xi\to\pm\infty}q_2(\xi,s,t)&\sim \at{\tfrac{1}{2}\ds^2 q^\pm_0}_{\sigma=0}\xi^2 
	+ \at{\ds q^\pm_1}_{\sigma=0} \xi
	+ \at{q^\pm_2}_{\sigma=0}.
	\end{align*}
Equipped with our general multiple scales matched asymptotic procedure for arbitrary geometries, we can now apply it to the volume-penalty method and derive its asymptotic behaviour in a straightforward step-by-step manner.

\subsection{Zeroth order: The tangential boundary layer}
\label{sec:zeroth}
\emph{Fluid problem} --- At leading order in $\eps$ we reproduce the fluid equations \cref{eq:ns-true}
	\begin{align}
    	\pd{u^+_0}{t} - \frac{1}{\Re} \lapof u^+_0 + \grad p^+_0 + u^+_0 \cdot \grad u^+_0 &= f, &     	
    	\divof{u^+_0} &= 0, &
    	u^+_0 &= 0 \quad \text{on} \quad \partial \Omega.
    \end{align} 
\emph{Solid problem} --- The leading order solid problem is 
	\begin{align}
    \lapof{p^-_0} &= \divof{f}, & u^-_0 &= 0.
    \end{align}
\emph{Boundary problem} --- The leading order boundary problem illustrates the structure of the subsequent orders.
The divergence condition reduces to purely the normal direction, and implies the normal velocity is constant. Matching to the solid problem requires this to be zero, giving a no-flux boundary condition for the fluid velocity,
    \begin{align}
    \dx u_{\sigma0} = 0 \qquad\implies\qquad u_{\sigma0} = 0 \qquad \implies \qquad u^+_{\sigma0}(0) = 0.
    \end{align}

{\textbf{Tangential velocity boundary layer}} --- 
We now analyse the tangential velocities with the leading order momentum equation, which takes the form
    \begin{align}\label{eq:inner-problem}
    (\Gamma - \dx^2) u_{\bot 0} &= 0.
    \end{align}
We emphasise this leading order operator, as it is key to understanding the volume-penalty method.
Noting the asymptotic behaviour of $\Gamma$, we can infer the asymptotic behaviour of the kernel.
The first solution $\V$ grows unphysically in the solid, and is constant in the fluid
    \begin{align}
    \V(\xi \to -\infty) &\sim a e^{-\xi}, & \V(\xi \to + \infty) &\sim 1.
    \end{align}
The \emph{physical} solution $\U$ decays exponentially into the solid, and grows linearly in the fluid
    \begin{align}
    \U(\xi \to -\infty) &\sim b e^{\xi}, & \U(\xi \to + \infty) &\sim \xi - \ell^*.
    \end{align}
{The shift $\ell^*$ is the \emph{displacement length} of the mask}.
Example solutions for different $\Gamma$ are shown in \cref{fig:inner-diagrams}.
The solid matching condition tells us that $u_{\bot 0}$ is some multiple of the $\U$ solution,
    \begin{align}
    u_{\bot 0} = A_{\bot 0} \, \U.
    \end{align}
The fluid matching condition instead requires zero gradient in the fluid.
This is only satisfied by the zero solution, implying no-{slip} boundary conditions at $\partial \Omega_s$ for the fluid velocities,
    \begin{align}
    u_{\bot 0}(\xi \to \infty) &\sim u^+_{\bot 0}(0) &
    \implies & &
    u_{\bot 0} &= A_{\bot 0} = u^+_0(0) = 0.
    \end{align}
The leading order problem is therefore equivalent to the incompressible Navier-Stokes equations with no-slip boundaries at $\sigma = 0$,
justifying convergence of the volume-penalty method.
We continue with the next order problem, the leading order error, and show how it can be controlled through judicious choice of $\Gamma$.

\subsection{First order: Displacement length error}
\label{sec:first}
\emph{Fluid problem} --- In the fluid we now have
	\begin{align}
    	\pd{u^+_1}{t} - \frac{1}{\Re} \lapof u^+_1 + \grad p^+_1 + u^+_0 \cdot \grad u^+_1 + u^+_1 \cdot \grad u^+_0 &= 0,&
    	\divof{u^+_1} &= 0,&
    	u^+_1 &= 0 \quad \text{on} \quad \partial \Omega.
    \end{align}
\emph{Solid problem} --- The solid then gives
	\begin{align}
    \lapof{p^-_1} &= 0, & u^-_1 &= 0.
    \end{align}
\emph{Boundary problem} --- The continuity equation for the boundary problem implies the first order  normal velocity is again constant, which is zero by matching with the solid boundary condition
    \begin{align}
    \dx u_{\sigma 1} - \Kbar \dx(\xi u_{\sigma0}) + \sgrad \cdot u_{\bot 0} &= \dx u_{\sigma 1} = 0 &
    \implies & &
    u_{\sigma1} &= 0.
    \end{align}
Hence the normal velocity is at most $\O(\eps^2)$ around the boundary.
The normal component of the momentum equation implies the leading order behaviour of the pressure is constant
    \begin{align}
    \frac{1}{\Re} \left(\Gamma u_{\sigma1} + \sdiv(\dx u_{\bot 0})\right) + \dx p_{0} + u_{\sigma0} \dx u_{\sigma0} &= \dx p_0 = 0 & &\implies & p_0 = p^\pm_0(0).
    \end{align}
    
{\textbf{Error for general masks} --- 
The tangential ($\bot$) momentum equation implies
    \begin{align}
    (\Gamma - \dx^2) u_{\bot 1} + \Kbar \dx u_{\bot 0} + \dx \sgrad u_{\sigma 0} + \Re u_{\sigma0} \dx u_{\bot 0} &= (\Gamma - \dx^2)\,u_{\bot 1} = 0.
    \end{align}
Matching to the solid problem implies $u_{\bot 1}$ is some multiple of $\U$
    \begin{align}
    u_{\bot 1} = A_{\bot 1}\, \U,
    \end{align}
and matching the fluid problem gives us the boundary condition for $u_{\bot 1}^+$,
    \begin{align}
    A_{\bot 1}\, \U &\sim u^+_{\bot1}(0) + \xi \ds u^+_{\bot0}(0) & &\implies & A_{\bot 1} &= \ds u^+_{\bot0}(0), &
    &\implies & u^+_{\bot1}(0) &= -\ds u^+_{\bot0}(0) \, \ell^*. \nonumber
    \end{align}
This shift $\ell^*$ is the {displacement length} of the mask},
and depends {only} on $\Gamma$.
Most masks give displacement length $\ell^*$ order unity (see \cref{sec:inner} for examples), implying the boundary condition $u^{+}_{\bot 1}(0) = \O(\eps)$.
This boundary perturbation of the fluid problem causes first order error $\O(\eps)$ {throughout} the fluid.

{\textbf{Error for optimal masks }} --- 
However, if we choose a mask with zero displacement length ($\ell^* = 0$), we obtain homogenous boundary conditions in the fluid,
    \begin{align}
    \left.u^+_1\right|_{\sigma=0} &= 0 \quad \text{on} \quad \partial \Omega_s, &
    \text{and} & &
    u^+_1 &= 0 \quad \text{on} \quad \partial \Omega.
    \end{align}
If the errors $u^+_1$ and $p^+_1$ are initially zero, the homogeneous linear momentum equation implies they will remain so for all time.
In reality, any small perturbation will grow over time (according to the Lyapunov time scale of the flow), but this is inherent to the chaotic nature of most fluid-flows, and not a limitation of this specific method.
The statistical behaviour of the fluid will be accurate to $\O(\eps^2)$ with this correction.
We show how to choose a mask with no displacement error in \cref{sec:inner}, but first analyse the second order problem, and in doing so find the emergence of two $\Re$ dependent asymptotic regimes.

\renewcommand{\W}[0]{\mathcal{W}}

\subsection{Second order: Different volume-penalty regimes}
\label{sec:second}
\emph{Fluid problem} --- We have an added inhomogeneity from the nonlinear first order convection term, which is zero for optimised masks,
	\begin{align}
	    \pd{u^+_2}{t} - \frac{1}{\Re} \lapof u^+_2 + \grad p^+_2 + \sum\nolimits_{j=0}^2 u^+_j \cdot \grad u^+_{2-j} &= 0, &
    	\divof{u^+_2} &= 0,&
    	u^+_2 &= 0 \quad \text{on} \quad \partial \Omega.
    \end{align}
\emph{Solid problem} --- We now have non-zero flow due to pressure gradients and forcing in the solid,
	\begin{align}
    \lapof{p^-_2} &= 0, & u^-_2 &= \Re \left(f -\grad p^-_0 \right).
    \end{align}
\emph{Boundary problem} --- 
The incompressibility equation simplifies since $u_{\sigma 1} = 0$, giving $u_{\sigma 2}$
    \begin{align}
    \dx u_{\sigma 2} + \sgrad \cdot u_{\bot 1} &= 0, &
    u_{\sigma 2} &= u_{\sigma 2}(-\infty) - \sdiv{(\ds u_0^+(0))} \int^\xi_{-\infty} \U \, d\zeta.
    \end{align}
We can simplify noting that at the boundary $\sigma = 0^\pm$,
    \begin{align}
    u_{\sigma 2}(-\infty) &= \Re (f_{\sigma} - \ds p^-_0)_{\sigma = 0^-}, &
    \sdiv{(\ds u_0^+(0))} &= \Re (f_{\sigma} - \ds p^+_0)_{\sigma = 0^+},
    \end{align}
where the first equality comes from matching to the solid solution, and the second equality comes from evaluating the normal component of the fluid momentum equation at $\sigma = 0$.
This behaves as a constant plus an exponential in the solid, and quadratically within the fluid,
    \begin{align}
    \text{as} \quad \xi \to -\infty \qquad u_{\sigma 2} &\sim \Re \left[ \left(f_{\sigma} (0^-) - \ds p^-_0(0) \right) (1 - b e^\xi)\right], \\
    \text{as} \quad \xi \to +\infty \qquad u_{\sigma 2} &\sim \Re \left[ \left(f_{\sigma} (0^+) - \ds p^+_0(0) \right) \frac{\xi^2}{2} + A_{\sigma 2} \xi + B_{\sigma 2}\right].
    \end{align}
The quadratic behaviour is consistent with the leading order problem,
and the linear behaviour matches onto zero for an optimised mask ($A_{\sigma 2} \propto \ell$).
However, the constant behaviour into the fluid depends on the jump in the pressure derivative (and external forcing) across the boundary, which depends on the details of the flow, and cannot be eliminated using a passive mask, so in general $B_{\sigma 2}$ is non-zero.
Hence, the boundary condition for the second order normal velocity in the fluid scales with the {Reynolds number} $\Re$ (for a given jump in the pressure).

We can also take the divergence of the momentum equation in the boundary region to obtain
    \begin{align}
    \dx^2 p_1 = \Kbar \dx p_{0} -\frac{1}{\Re}\dx \Gamma u_{\sigma 2} + \dx f_\sigma.
    \end{align}
Integration shows a discontinuity in the normal pressure gradient across the mask boundary.

We can then determine the tangential velocity using the tangential momentum equation,
    \begin{align}\label{eq:second-order}
    (\Gamma - \dx^2) u_{\bot 2} = - \Kbar \dx u_{\bot 1}  + \Re (f_\bot - \sgrad p_{0}) = \mathcal{R}.
    \end{align}
This is solved using variation of parameters,
    \begin{align}
    u_{\bot 2} &= (\mathcal{Q}+c) \, \U, &
    \mathcal{Q} &\equiv \int -\frac{\int_{-\infty}^\eta \mathcal{R}{\U} d\zeta}{\U^2}d\eta, 
	\end{align}
where $\mathcal{R}$ is the right hand side of \cref{eq:second-order} and $c$ is a constant of integration.

\begin{figure}[ht]
        \centering
        \includegraphics[width=.75\linewidth]{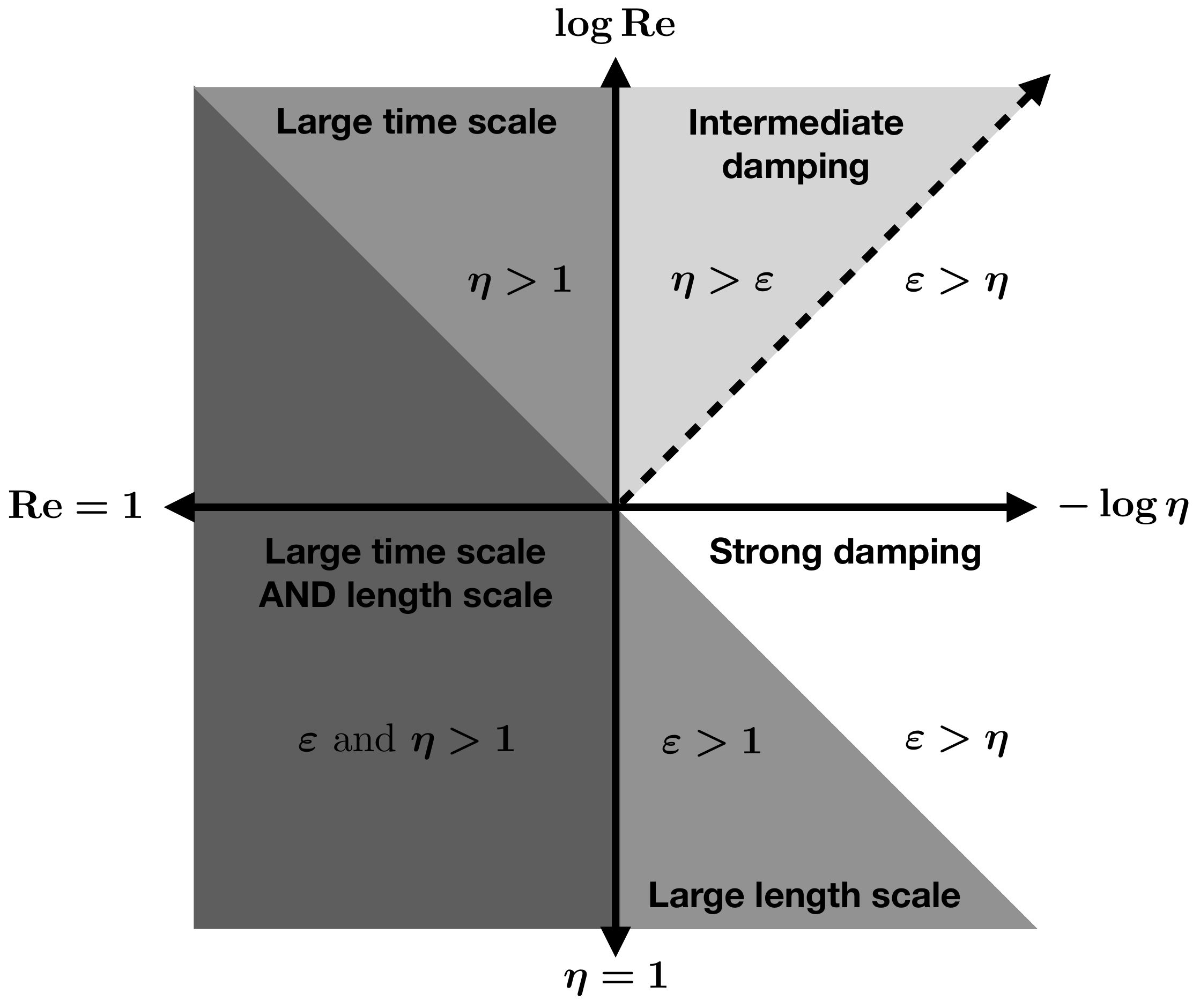}
        \caption{Different volume-penalty regimes shown on a log-log plot  of the Reynolds number $\Re$ (vertical axis) and damping {rate} $\eta^{-1}$ (horizontal axis). The regions are distinguished by the relative size of the damping length scale $\eps$ and damping time scale $\eta$. Increased damping (smaller $\eta$) corresponds to moving right on the diagram. The darkest regions are where at least one of $\eta$ or $\eps$ are greater than 1, and the solution is incorrect. Our correction is useful in the strong damping regime (white), which is the limiting regime as $\eta \to 0$ for any fixed $\Re$. For $\Re > 1$, an intermediate regime can arise (light gray), in which the damping time scale dominates the damping length scale $1 > \eta > \eps > \Re^{-1}$.
          }
        \label{fig:regime-diagrams}
    \end{figure}

{\textbf{Reynolds dependent error regimes}} --- 
The second order error in the fluid scales with the Reynolds number $\Re$, and so is proportional to $\eta$, while the first order error scales with $\eps$.
For sufficiently large $\Re > \eps^{-1}$, the asymptotic separation is invalid, and the ``second order'' time scale error $\O(\Re \, \eps^2) = \O(\eta)$ dominates the ``first order'' {length} scale error $\O(\eps)$.
This leads to several asymptotic damping regimes.
We summarise these regimes in terms of $\log\eta$ and $\log\Re$ in \cref{fig:regime-diagrams}.

If $\eta > 1$ (the left half plane), the fluid is damped slowly and the velocity in the solid is large.
Likewise if $\eps > 1$ (the lower left diagonal half-plane), damping cannot overcome the fluid viscosity.
It is possible to have either $\eps$ or $\eta$ be small and the other large, but \emph{both} must be small to approximate a solid.

Ignoring these large error regimes, we can have \emph{intermediate} damping  (light gray) when the Reynolds number is large enough to promote the second order system to dominance ($1 > \eta > \eps > \Re^{-1}$).
This regime gives the appearance of $\O(\eps^2)$ convergence even for unoptimised masks, but this is only {temporary}, as seen in \cite{AngotPenalizationMethodTake1999}.

Once the damping is sufficiently large ($1 > \Re^{-1} > \eps > \eta$) the damping length scale $\eps$ always dominates the time scale $\eta$.
This is the \emph{strong} damping regime (white), for which it is {necessary} to apply mask corrections to promote the convergence to $\O(\eps^2)$.

{\textbf{Reynolds dependent cost regimes}} --- 
We can analyse the computational cost regimes for large $\Re$ by comparing the damping length scale $\eps$ and time scale $\eta$ with the intrinsic Kolmogorov length scale $\lambda$ and CFL time scale $t$ in three dimensional turbulence,
    \begin{align}
    \Re &= \frac{UL}{\nu}, & 
    \eta &= \frac{\tau U}{L}, & 
    \eps &= \sqrt{\frac{\eta}{\Re}}, &  
    t &= \Re^{-3/4}, & 
    \lambda &= \Re^{-3/4}.
    \end{align}
The grid resolution $\Delta x$ and time step size $\Delta t$ must resolve the smallest scale,
    \begin{align}
    \Delta x &= \min (\eps, \lambda) & \text{and}& & \Delta t &= \min(\eta,t).
    \end{align}
The computational cost $C$ is proportional to the degrees of freedom $\Delta x^{-3} \Delta t^{-1}$.
Scaling the damping in terms of the Reynolds number $\eta = \Re^{-\alpha}$, we find several cost regimes for $\Re \gg 1$,
    \begin{align}
    \alpha &< 1/2       & &\implies & \lambda &= t < \eps < \eta & &\implies & C &\sim \, \lambda^{-3}t^{-1} = \Re^{3},\\
    1/2 < \alpha &<3/4  & &\implies & \eps &< \lambda = t < \eta & &\implies & C &\sim \, \eps^{-3}t^{-1} = \eta^{-3/2}\Re^{3/4}, \\
    3/4 < \alpha &< 1   & &\implies & \eps &< \eta < \lambda = t & &\implies & C &\sim \, \eps^{-3}\eta^{-1} = \eta^{-5/2}\Re^{3/2},\\
    1 < \alpha &        & &\implies & \eta &< \eps < \lambda = t & &\implies & C &\sim \, \eps^{-3}\eta^{-1} = \eta^{-5/2}\Re^{3/2}.
    \end{align}
Volume penalised boundaries can thus be added at no extra cost with an error of $\Re^{-1/2}$.
Beyond this regime $\alpha > 1/2$, the computational cost scales more aggressively in $\eta$.
This has given the impression that low $\Re$ simulations are difficult using volume penalisation.
However, smaller $\Re$ makes it easier to achieve strong damping where mask corrections can now allow $\O(\eps^2)$ convergence.
We therefore extend the practicability of the volume-penalty method to low $\Re$.
We finally note that this analysis is premised on being in the strong damping regime, where the second order error is smaller than the first order error.

\subsection{Force and torque calculations}
We now relate the integral of the penalty term to the force and torque calculations.
For a general mask, this integral can be rewritten using the penalised Navier-Stokes equations
\begin{align}
	\int_\Omega \frac{1}{\Re\eps^2} \Gamma u \, dV &= \int_\Omega \divof{\Sigma} - \frac{D u}{D t} + f \, dV, &
	\Sigma &\equiv - p I + \frac{1}{\Re}\left(\nabla u + \nabla u^\top\right).
    \end{align}
We use a body-centered coordinate frame and account for the acceleration terms in the body forcing $f$.
For a mask with support on $\Omega\setminus\Omega^+$, which is $\O(\eps)$ larger than $\Omega_s$, we can restrict the integral to this region, with exponentially small error for non-compact masks like $\Gamma_{\tanh}$.
    \begin{align}
    \int_\Omega \frac{1}{\Re\eps^2} \Gamma u \, dV = \int_{\Omega\setminus \Omega^+} \divof{\Sigma} - \frac{D u}{D t} + f \, dV.
    \end{align}
The divergence theorem relates this to the stress along $\partial \Omega^+$, plus corrections
    \begin{align}
    \int_\Omega \frac{1}{\Re\eps^2} \Gamma u \, dV = \int_{\partial \Omega^+} \hspace{-1em} \Sigma \cdot \nhat \,dS - \int_{\partial \Omega_{in}} \hspace{-1em} \Sigma \cdot \nhat \,dS + \int_{\Omega^+} - \frac{D u}{D t} + f \, dV.
    \end{align}
We can compare the penalised stress $\Sigma$ with the true stress $\Sigma_0$
    \begin{align}
    \int_\Omega \frac{1}{\Re\eps^2} \Gamma u \, dV = \int_{\partial \Omega^+} \hspace{-1em} (\Sigma_0 + \Sigma - \Sigma_0) \cdot \nhat \,dS - \int_{\partial \Omega_{in}} \hspace{-1em}\Sigma \cdot \nhat \,dS + \int_{\Omega^+} - \frac{D u}{D t} + f \, dV.
    \end{align}
We can then relate the true stress on the mask support boundary $\partial \Omega^+$ to the true stress on the true boundary $\partial \Omega_s$ by Stokes' theorem with the unpenalised equation for $u_0, p_0$
    \begin{align}
	\int_{\partial \Omega_s} \Sigma_0 \cdot \nhat \, dS &= \int_\Omega \frac{1}{\Re\eps^2} \Gamma u \, dV + \int_{\partial \Omega_{in}} \hspace{-1em}\Sigma \cdot \nhat \,dS - \int_{\Omega_s} f \, dV - \int_{\partial \Omega^+} (\Sigma - \Sigma_0) \cdot \nhat \,dS \\
	 & \quad - \int_{\Omega^+ \setminus\Omega_s} \frac{D u_0}{Dt} \, dV + \int_{\Omega^+} \frac{D u}{D t} \, dV .\nonumber
    \end{align}
Both the stress on interior no-slip boundaries (as in \cref{sec:unsteady-cylinder}) and the $f$ volume integral are at order $\eps^0$, so must be subtracted off (note that $f$ is integrated over exactly $\Omega_s$).
The leading stress error $\Sigma-\Sigma_0$ is $\O(\eps)$ at $\Omega^+$ for general masks, as both $p^+_1$ and $u^+_1$ are non-zero in the fluid, implying $\O(\eps)$ convergence of the force calculation.

However, an optimal mask ensures that $p^+_1$ and $u^+_1$, and hence $\Sigma - \Sigma_0$ are accurate to $\O(\eps^2)$ at the mask boundary $\partial \Omega^+$, so this term is now $\O(\eps^2)$ accurate.
The remaining acceleration terms are $\O(\eps^2)$ in general as both $u_0$ and $\eps u_1$ are $\O(\eps)$ in a region $\Omega^+ - \Omega_s$ with  support $\O(\eps)$.
The $\O(\eps)$ error $u_1$ decays exponentially fast into the solid, so also has support and magnitude of $\O(\eps)$, integrating to $\O(\eps^2)$.
The $\O(\eps^2)$ error $u^-_2$ and $p^-_2$ have support throughout the solid region, but again integrate to only $\O(\eps^2)$.
Hence, by choosing an optimal mask, the dominant stress error can be eliminated and the force calculation accuracy promoted to $O(\eps^2)$.

The torque calculation proceeds similarly.
We recover the correct stress at $\partial \Omega^+$ to $\O(\eps^2)$, and we need to show this is only $\O(\eps^2)$ different from the true stress at $\partial \Omega_s$.
This follows from the divergence theorem applied to $r \times (\divof \Sigma_0)$, which relies on the fact that the stress tensor is symmetric, and that the gradient of the coordinate vector is the identity tensor.
Following the same procedure as for the force calculation, one finds $\O(\eps)$ accuracy for the torque for general masks, and $\O(\eps^2)$ accuracy for optimal masks.

\section{Testing optimal masks}
In the previous section, we derived the necessary conditions for optimal masks.
We now investigate these conditions, and test the performance of our corrections in a series of increasingly general test problems.

\subsection{Calculating optimal masks}
\label{sec:inner}
\begin{figure}[h]
        \centering
        \includegraphics[width=.9\linewidth]{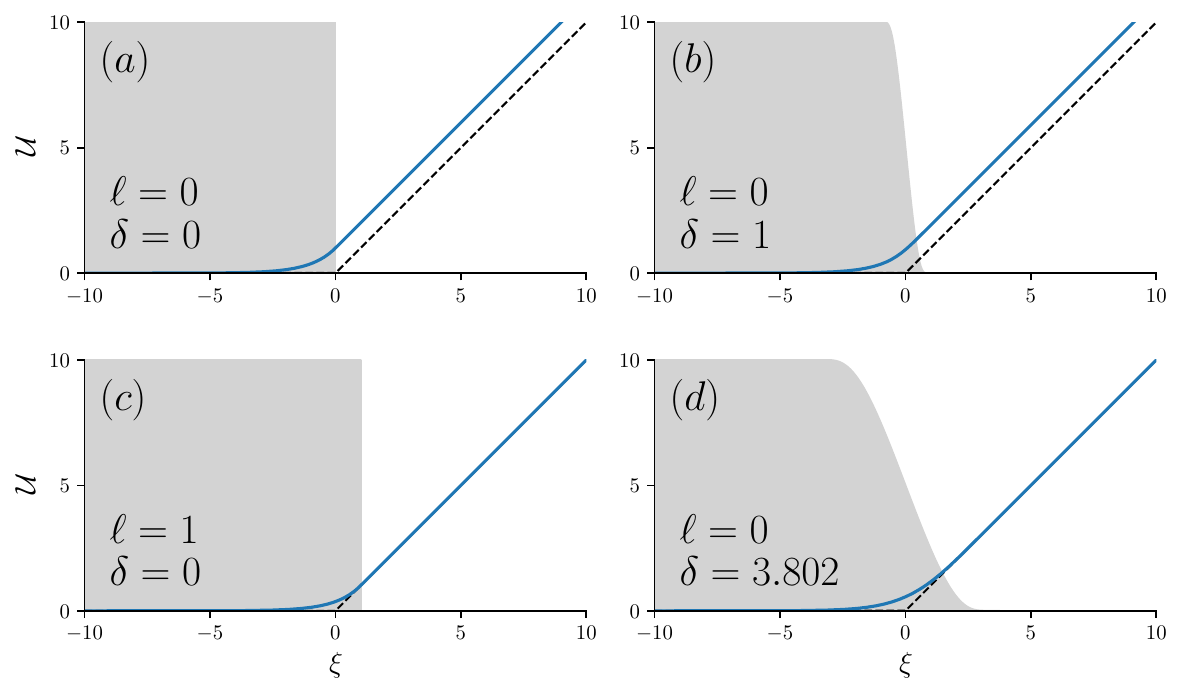}
        \caption{
          Different mask choices (grey) give different velocity profiles (blue) compared to the ideal solution (dashed black) in the boundary layer problem \cref{eq:inner-problem}.
          Figure $(a)$ shows the conventional discontinuous mask, 
          with zero smoothing or shifting, results in a displacement length $\ell^* = -1$.
          Plot $(b)$ shows an unoptimised compact erf profile using smoothness $\delta = 1$ gives $\O(1)$ displacement length.
          Plot $(c)$ optimally shifts the standard mask by $\ell = 1$ to eliminate the displacement length error.
          Plot $(d)$ shows optimal smoothing $\delta^*$ (\cref{tab:inner-smoothings}) can also eliminate the displacement length error.
          }
        \label{fig:inner-diagrams}
    \end{figure}

We first examine the construction of optimal masks for the inner problem of \cref{sec:asymptotics} \cref{eq:inner-problem},
    \begin{align*}
    \dx^2 \U(\xi) - \Gamma(\xi)\, \U(\xi) &= 0, & \U(\xi\to -\infty)&\to 0, & \U(\xi\to +\infty) &\to 1.
    \end{align*}
We plot several solutions $\U(\xi)$ to various mask function choices $\Gamma$ in \cref{fig:inner-diagrams}.
    
\subsubsection{Conventional volume penalisation}
Conventional masks in \cref{eq:inner-problem} lead to non-zero displacement length,
	\begin{align}\label{eq:inner-mask-disc}
	\Gamma_{0,0}(\xi) &= \begin{cases}
 		1 & \xi < 0,\\
 		0 & \xi > 0,
 	\end{cases}
	&
	&\implies
	&
	\U(\xi) &= 
	\begin{cases}
	e^{\xi} & \xi < 0,\\
	1 + \xi & \xi > 0.
	\end{cases}
	\end{align}
From \cref{sec:asymptotics} the non-zero displacement length $\ell^* = -1$ induces an $\O(\eps)$ error in the fluid solution, agreeing with previous estimates of the accuracy of the volume-penalty method.
We now examine mask corrections to eliminate this displacement and improve convergence.

\subsubsection{Shifting the mask}
The simplest correction is to shift the mask distance $\ell = 1$ into the fluid (\cref{fig:inner-diagrams} ($a$) vs \cref{fig:inner-diagrams} ($b$)).
	\begin{align}\label{eq:inner-mask-disc-shift}
	\Gamma_{\ell^*,0}(\xi) &= \begin{cases}
 		1 & \xi < 1,\\
 		0 & \xi > 1,
 	\end{cases} &
 	&\implies &
	\U(\xi) = 
	\begin{cases}
	e^{\xi - 1} &\xi < 1,\\
	\xi & \xi > 1.
	\end{cases}
	\end{align}	
This simple refinement eliminates the displacement error, and \cref{sec:asymptotics} shows this improves the global accuracy $\E_1$ and drag and torque errors $\Delta F$ and $\Delta T$ to $\O(\eps^2)$.

\subsubsection{Smoothing the mask}
We now demonstrate a key finding of the paper: there exist smooth masks that eliminate the displacement length error.
Late in the review stage we became aware of work by Poulsen and Voorhes which independently performed a similar analysis \cite{PoulsenSmoothedBoundaryMethod2017}.
We consider a $\tanh$ type mask,
    \begin{align}
    \Gamma(\xi) = \frac{1}{2}\left(1 - \tanh\left(\frac{2\xi}{\delta}\right)\right).
    \end{align}
We define transformed coordinates and variables and a rescaled smoothness $n$
    \begin{align}
    z(\xi) &\equiv \frac{1+\tanh\left(\frac{2\xi}{\delta}\right)}{2},
     & \U(\xi) &= U(z(\xi)),
     & n \equiv \frac{\delta}{4},
    \end{align}
for which the problem becomes
    \begin{align}
    (1 - z) z^2 U''(z) + (1 - 2 z) z U'(z) - n^2 U(z) = 0.
    \end{align}
The solution which satisfies $U(0) = 0$ can be written as a Frobenius series
    \begin{align}
    U(z) =  n\sum_{k=0}^{\infty} \frac{\Gamma (k+n) \Gamma (k+n+1) }{\Gamma
   (k+2 n+1)} \frac{z^{k+n}}{k!} = n \mathrm{B}(n,n+1) \, z^n {}_2F_1(n, n+1, 2n+1, z),
    \end{align}
where here $\Gamma(n) = (n-1)!$ represents the actual gamma function, $\mathrm{B}$ is the beta function, and ${}_2 F_1$ is the hypergeometric function.
Considering the behaviour of the function as $z\to 1$ from below, and applying the transformation, we can show that as $\xi \to \infty$
    \begin{align}
    \U(\xi) &\sim \xi - \ell(n),
	&
    \ell(n) &\equiv 1 + 2n(\psi(n) + \gamma),
    &
    \psi(n) &\equiv - \gamma + \sum_{k=0}^{\infty} \frac{n-1}{(k+n)(k+1)},
    \end{align}
where $\gamma$ is the Euler--Mascheroni constant and $\psi(n)$ is the digamma function $\Gamma'(n)/\Gamma(n)$.
We can calculate the value of $n$ (and hence $\delta$) for which the displacement length $\ell(n)$ is zero
    \begin{align}
    \ell(n) = 0, \quad \mathrm{for} \quad n \approx 0.662057 \quad \text{or} \quad \delta \approx 2.64822828 \, . 
    \end{align}
This is the optimal {smoothing} for an unshifted tanh mask which cancels the displacement error.

We now extend this work to a simple calculation for any mask profile.

\subsubsection{General corrections via Riccati transform}

    \begin{figure}[ht]
        \centering
        \includegraphics[width=.6\linewidth]{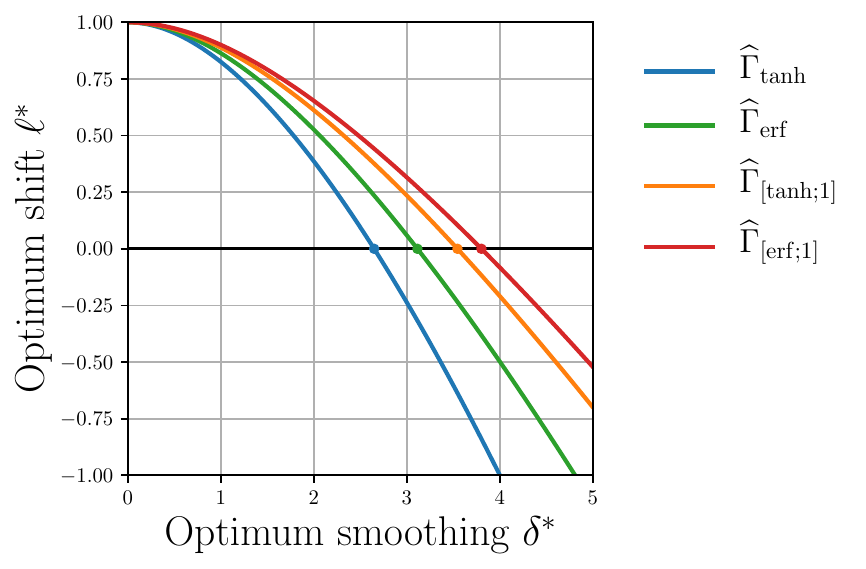}
        \caption{The optimum shift $\ell^*$ as a function of smoothing $\delta^*$ for four families of mask function. The limit $\delta \to 0$ reproduces the optimum discontinuous shift $\ell^* = 1$, and the zero shift smoothing occurs at the $\ell=0$  curve intercepts, given in \cref{tab:inner-smoothings}.}
        \label{fig:inner-shifting-vs-smoothing}
    \end{figure}

\newcommand\Tstrut{\rule{0pt}{2.6ex}}         
\newcommand\Bstrut{\rule[-0.9ex]{0pt}{0pt}} 

    \begin{table}[ht]
	\caption{Optimal zero-shift smoothing $\delta^*(\ell=0)$ for different normalised mask profiles $\Ghat$.}
	\begin{center}
	\begin{tabular}{l|cccc}
	{Mask profile} & $\Ghat_{\erf}$ & $\Ghat_{[\erf;1]}$ & $\Ghat_{\tanh}$ & $\Ghat_{[\tanh;1]}$\\
	{Optimum smoothing} & 3.11346786 & 3.80171928 & 2.64822828 & 3.54403048
	\end{tabular}
	\end{center}
	\label{tab:inner-smoothings}
    \end{table}%

We now calculate optimal smoothing and shifting for general mask functions.
In current form this requires solving second order boundary value problems.
We greatly simplify by using the shifted and scaled coordinate $\zeta$.
This transforms the velocity $\U(\xi)$ to the velocity $u(\zeta)$, and a general mask $\Gamma_{\ell,\delta}(\xi)$ to a normalised mask $\Ghat(\zeta)$.
We then apply a Riccati transform $R$ to the rescaled \cref{eq:inner-problem} to derive an initial value problem for $R(\zeta)$,
	\begin{align}
	\zeta &= \frac{\xi - \ell}{\delta}, &
    R(\zeta) &= \frac{\dot{u}(\zeta)}{u(\zeta)} = \frac{\delta \, \U'(\xi)}{\U(\xi)}, &
	R'(\zeta) + R^2(\zeta) &={\delta}^2 \Ghat(\zeta).
	\end{align}
For a mask equal to one for $\zeta < -c$ we have an exponential solution $\U(\xi) \propto e^\xi$.
Differentiating gives an initial condition for the Riccati equation at $\zeta = c$,
	\begin{align}
	R(\zeta=-c) &= {\delta} \frac{\U'(\xi)}{\U(\xi)} = {\delta}.
	\end{align}
Solving forward to $R(\zeta=c)$, we then enforce the optimal linear profile $u(\zeta) = \zeta$ beyond $\zeta > c$ to determine the ideal shift $\ell^*$ in terms of smoothing $\delta$,
    \begin{align}
    u(\zeta=c) &= \zeta, &
    &\implies &
    R(\zeta=c) &= \frac{\delta}{\ell^* + c \delta}, &
    &\implies & 
	\ell^*(\delta) &= \br{\frac{1}{R(c)} - c}\delta.
	\end{align}

This transformed system is much easier to solve than the original second order boundary value problem.
We solve it using simple Runge-Kutta integrators of the scipy python library.
We can then extend to noncompact smooth functions by choosing $c$ sufficiently large.

To illustrate the optimal parameters we plot the (scaled) shift $\ell^*$ as a function of smoothing $\delta$ for four normalised masks in \cref{fig:inner-shifting-vs-smoothing}.
The discontinuous limit $\delta \to 0$ reproduces the optimum discontinuous shift $\ell^* = 1$, and as the smoothing increases, the optimum shift decreases.

An important case is the `zero-shift' optimum smoothing.
These are the $\delta^*$ intercepts in \cref{fig:inner-shifting-vs-smoothing}, given in \cref{tab:inner-smoothings}.
This choice is numerically motivated by balancing requirements of resolving the mask and the velocity field.
Moving away from the optimum curve, the mask is either too `large' ($\ell > \ell^*(\delta)$), or too small ($\ell < \ell^*(\delta)$).

\subsection{Couette flow: Optimal parameters to minimise total error}
    \begin{figure}[ht]
    \centering
    \includegraphics[width=.9\linewidth]{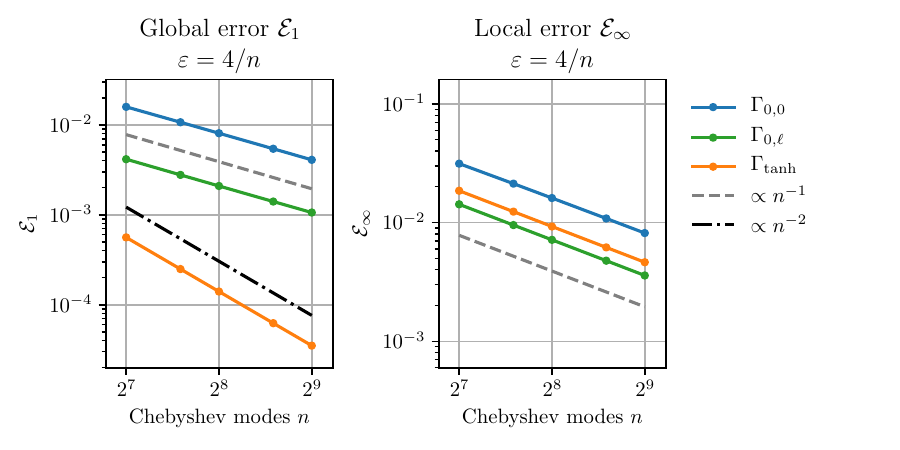}
    \caption{We plot the global $\E_1$ and local $\E_\infty$ error of the solution to Couette flow as a function of Chebyshev modes $n$, for the choice $\eps=4/n$. 
    We compare three strategies for optimisation: standard discontinuous masks $\Gamma_{0,0}$, optimal shifted discontinuous masks $\Gamma_{0,\ell}$, and optimal smoothed tanh masks $\Gamma_{\tanh}$.
    Optimal smoothed masks achieve second order numerical convergence of the total error, in contrast to discontinuous masks.
    }
    \label{fig:couette-numerical-error}
    \end{figure}
We have outlined the two sources of error of the volume-penalty method, mathematical \emph{model} error, and numerical \emph{discretisation} error.
The question faced by practitioners is, given some fixed computational resources, how can one minimise the \emph{total} error?
Damping influences the model error and discretisation error in contrasting ways.
Weak damping means the body does not approximate a solid, worsening model error.
Strong damping induces the small length scales that worsen the discretisation error if not resolved.
There must exist some optimum between these extremes which minimises the sum of the model and discretisation error.

We examine the numerical performance of our corrections in a specific test problem: a Chebyshev discretisation of Couette flow.
The mathematical problem is a rescaling of the inner problem which solves for the tangential velocity $v$ near a volume-penalised boundary with a prescribed unit velocity gradient into the fluid 
	\begin{align}
	v''-\frac{\Gamma}{\eps^2} v &= 0, & v'(-1) &\to 0, & v'(1) = 1.
	\end{align}
	
We solve this problem numerically using Chebyshev polynomials in Dedalus \cite{BurnsDedalusFlexibleFramework2020} for five choices of numerical resolution, $n = 128,192,256,384,$ and $512$ Chebyshev modes.
For each resolution $n$ we compare three mask profile strategies: a standard discontinuous mask, an optimal shifted discontinuous mask $\ell = \eps$, and an optimal unshifted smoothed mask $\delta = \delta^* \eps$.
We choose $\eps = 4/n$ to balance the model and discretisation errors.
We then measure the total error as the $L^1$ and $L^\infty$ norm of the difference between the numerical solution and the ideal solution $v_0(x) = x$ in the fluid domain, on a resolved grid of 16 times higher resolution than the grid used to calculate the solution.

We plot the results in \cref{fig:couette-numerical-error}, which shows that smoothed masks achieve lower total error (in $L^1$ norm) than any discontinuous mask, for a given numerical resolution.
Shifted masks achieve slightly better $L^\infty$ error, but at the cost of large Gibbs oscillations that reduce numerical convergence to first order in the fluid domain.
We caution that the exact errors are dependent on different numerical discretisations.
However, the qualitative findings are general:
optimal smoothed masks improve the volume penalty method to second order numerical convergence.

\subsection{Poiseuille flow: Body forces with volume penalisation}
\begin{figure}[ht]
        \centering
        \includegraphics[width=.99\linewidth]{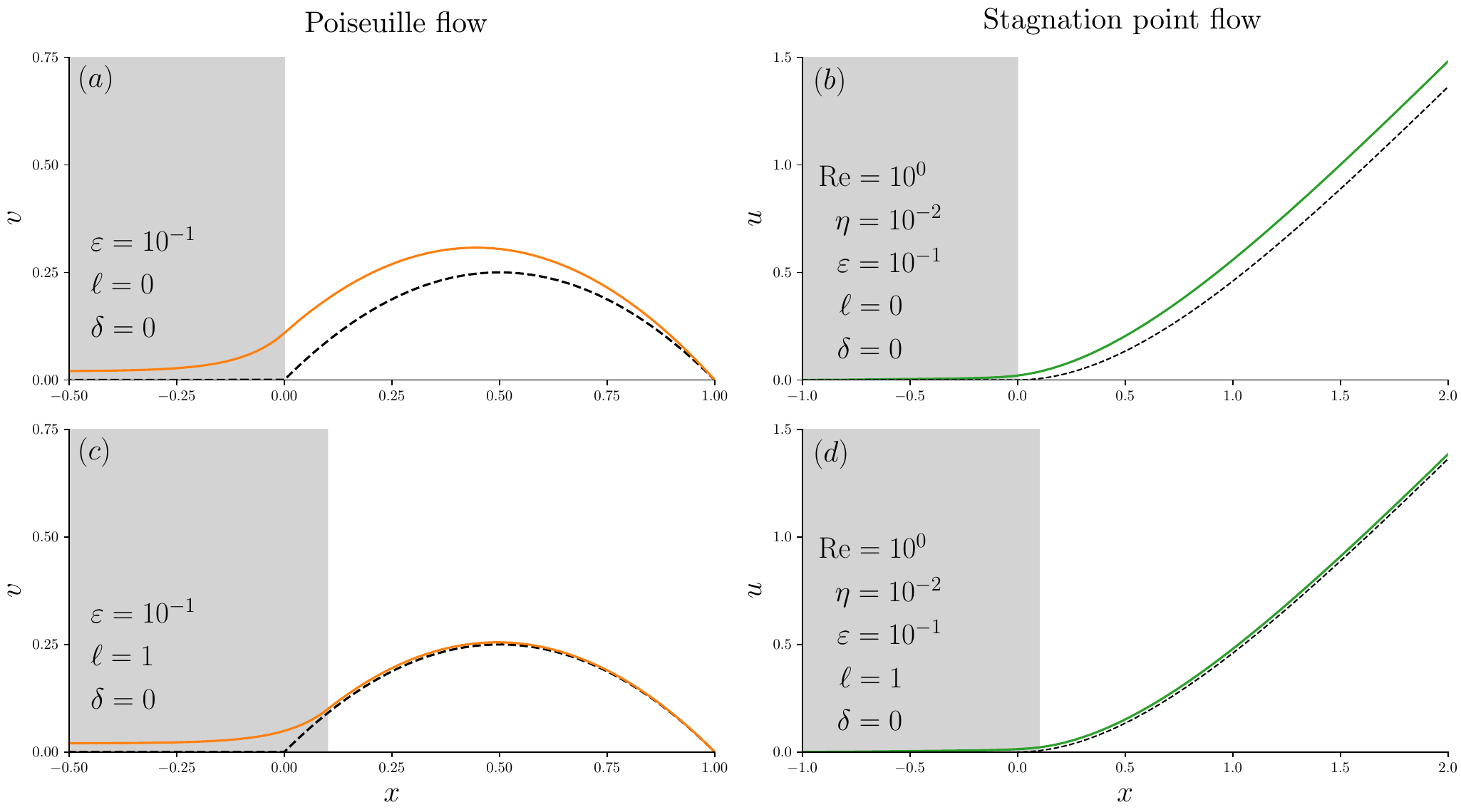}
        \caption{
        Figure $(a)$ plots the wall-tangential velocity $v$ (orange) of Poiseuille flow between a true no-slip boundary (at $x=1$), and a infinite volume penalised boundary beyond $x<0$, for damping length $\eps = 10^{-1}$.
        Figure $(b)$ plots the wall-normal velocity $u$ (green) of viscous stagnation point flow near a finite volume penalised boundary for $-1<x<0$ at $\Re = 1$, $\eta=10^{-2}$, which gives damping length $\eps = 10^{-1}>\eta$, corresponding to strong damping.
        Both plots use the standard unoptimised discontinuous mask $\ell = \delta = 0$ and give $\O(\eps)$ displacement error with respect to the ideal solutions (dashed black).
        Figure $(c)$ and $(d)$ apply zero-displacement length masks to improve the error to $\O(\eps^2)$ in the fluid.
          }
        \label{fig:poiseuille-stagnation-diagram}
    \end{figure}

We now test the robustness of optimal masks to forcing terms in Poiseuille flow, steady pressure driven flow between two no-slip walls.
This solution tests the robustness of optimal masks to forcing terms in the problem.
The penalised problem reduces to a single equation for the tangential velocity $v$ in terms of the wall-normal direction $x$,
	\begin{align}
	v''-\frac{\Gamma}{\eps^2} v &= -2, & v'(x\to-\infty) &\to 0, & v(1) = 0.
	\end{align}
A true wall replaces damping with no-slip boundary conditions at $x=0$.
In the fluid, the flow behaves parabolically,
and there is an $\O(\eps^2)$ flow due to pressure within the solid (\cref{fig:poiseuille-stagnation-diagram}).
However, the same $\O(\eps)$ displacement length occurs for conventional masks.
Optimised masks eliminate the $\O(\eps)$ displacement length, promoting the average error $\E_1$ to $\O(\eps^2)$, though the maximum local error $\E_\infty$ at $x=0$ remains $\O(\eps)$ (\cref{fig:poiseuille-stagnation-diagram} ($c$)).
We note that though it is possible to further reduce the fluid error outside of the mask, it is not possible to reduce the global error below $\O(\eps^2)$, or the local error below $\O(\eps)$.
This is because the error in the boundary region begins to dominate the error in the fluid region.

\subsection{Viscous stagnation-point flow: Reynolds dependent damping}
Flow about a stagnation point adds two important physical effects;
flow normal {and} tangent to a boundary, and equation nonlinearities parameterised by the Reynolds number $\Re$.
We can use symmetry to reduce the problem to a single third order nonlinear differential equation for the wall-normal velocity $u$ in the wall-normal coordinate $x$.
Adding a finite volume penalisation boundary, we now solve for the penalised normal velocity $u$,
	\begin{align}
	u'^2 - u'' u + \frac{1}{\Re}u''' -\frac{\Gamma}{\eta} u' &= 1, &
	u(-1) &= 0, &
	u'(-1) &= 0,&
	u'(\infty) &= 1.
	\end{align}
We cannot extend our volume penalised region indefinitely, due to the linear pressure driven flow in this region, which causes the velocity at $x=0$ to increase as the true boundary is moved further away.
We instead enforce no slip boundary conditions on $u$ and $u'$ at the non-dimensional location $x=-1$.
The true profile replaces damping with no-slip conditions applied at $x=0$.
We solve the system numerically using Dedalus \cite{BurnsDedalusFlexibleFramework2020}.

This simple system reproduces the distinction between normal and tangential velocity errors noted in \cref{sec:asymptotics}.
The combination of no-slip boundaries and incompressibility imply zero {gradient} in the normal velocity at the boundary  \cref{fig:poiseuille-stagnation-diagram} ($b$), ($d$).
Hence the normal velocity is always $\O(\eps^2)$ in the boundary region, whereas the tangential velocity can be $\O(\eps)$.

\begin{figure}
        \centering
        \includegraphics[width=.9\linewidth]{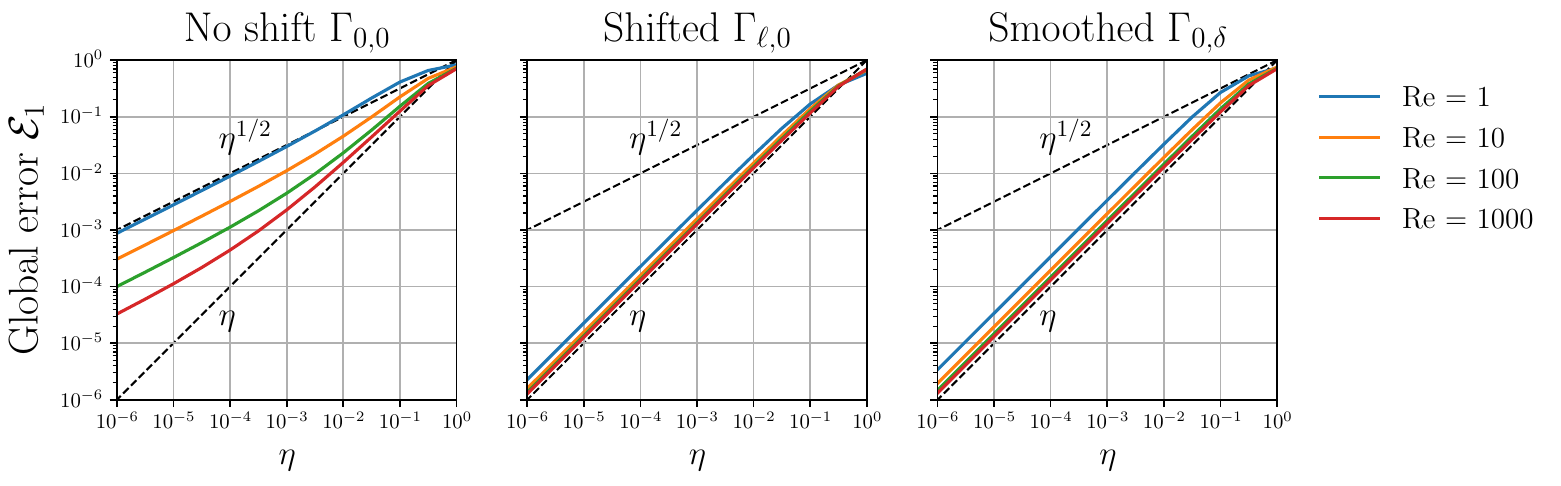}
        \caption{Log-log plots of global error $\E_1$ as a function of $\eta = \Re \eps^2$, for different choices of $\Re = 1,10,10^2,10^3$ in volume penalised stagnation point flow.
        Three choices of mask function are shown, being unshifted discontinuous mask $\Gamma_{0,0}$, optimum shifted discontinuous mask $\Gamma_{\eps,0}$, and optimum smoothed mask $\Gamma_{0,\delta}$ for the compact error function mask.
        A clear transition between intermediate and strong damping occurs for the unshifted discontinuous mask.
        The local error behaves similarly for stagnation point flow.
        }
        \label{fig:stagnation-Re-E1}
    \end{figure}

Stagnation point flow also contains nonlinearity and the Reynolds number $\Re$, leading to an additional damping regime.
As in \cref{sec:asymptotics}, we see a transition from intermediate to strong damping past $\eta < \Re^{-1}$ in a plot of the global error $\E_1$ of $u$ as a function of $\eta$ for different $\Re$ in \cref{fig:stagnation-Re-E1}.
While unoptimised masks can achieve temporary $\O(\eta)$ accuracy for intermediate damping, only zero displacement length masks can achieve $\O(\eta)$ convergence in large damping regimes (\cref{fig:poiseuille-stagnation-diagram} ($d$)).
This peculiar, if fortunate, high Reynolds number behaviour was observed experimentally by Angot et al. \cite{AngotPenalizationMethodTake1999}, contradicting the theoretical convergence predicted by Angot \cite{AngotAnalysisSingularPerturbations1999} ($\O(\eps^{1/2})$) and Carbou and Fabrie \cite{CarbouBoundaryLayerPenalization2003} ($\O(\eps)$).

\section{Flow past a rotating cylinder: Moving objects in 2D}\label{sec:unsteady-cylinder}

We now examine 2D incompressible flow past a rotating cylinder at moderate Reynolds number.
This system includes damping, viscosity, non-linearity, curvilinear geometries, and unsteady linear and rotational acceleration of the flow and body respectively, leading to non-zero unsteady horizontal drag $F_x$, vertical drag $F_y$, and torque $T$.
We simulate both a true no-slip cylinder and various volume penalised approximations using the open-source spectral PDE solver Dedalus \cite{BurnsDedalusFlexibleFramework2020} to determine the local and global error of each field as well as the physically relevant drag and torque accuracy.
We compare the performance of standard unoptimised masks with optimally shifted and optimally smoothed masks, validate our prescriptions, and apply a simple Richardson extrapolation scheme to further accelerate the convergence of optimised masks to $\O(\eps^4) \sim \O(\eta^2)$ accuracy.
The configurations are tested with geometry conforming spectral elements at {$\Re = 200$}, to examine the mathematical convergence in isolation of numerical error.
First we examine the mathematical problem.

\subsection{Mathematical formulation}
We simulate the volume penalised cylinder of unit radius on an annular domain $(\theta,r) \in [0,2\pi]\times[R_1,R_2]$. 
We do not include the center to avert the coordinate singularity.
We then solve the initial value problem with zero initial conditions, and governing equations
    \begin{align}
    \pd{(ru)}{r} + \pd{v}{\theta} &= 0,\\
    q - \pd{(rv)}{r} + \pd{u}{\theta} &= 0,\\
    r^2 \pd{u}{t} + \frac{1}{\Re} \pd{q}{\theta} + r^2\pd{p}{r} &= + rvq - \frac{1}{\eta} r^2 \Gamma(r) u,\\
    r^2 \pd{v}{t} - \frac{1}{\Re} (r\pd{q}{r} - q) + r \pd{p}{\theta} &= - ruq - \frac{1}{\eta} r^2 \Gamma(r) (v - r\Omega(t)).
    \end{align}
This is a slight reformulation of standard incompressible hydrodynamics in polar coordinates.
We solve for the radial and polar velocity $u$ and $v$, and have introduced the variables $p$ and $q$ to write the problem in first order form (as required for Dedalus), which are related to the true pressure $P$ and vorticity $\omega_z$ by
    \begin{align}
    p &= P + \frac{1}{2}(u^2 + v^2), & q &= r\omega_z.
    \end{align}  
We also multiply out factors of $r^{-1}$ to improve the convergence of spectral discretisations of the equations.
The outer boundary at $r = R_2$ applies uniform flow boundary conditions accelerating from rest
    \begin{align}
    u(R_2) &= g(t) \cos \theta, & v(R_2) &= -g(t) \sin \theta, & \text{where } g(t) = \Gamma_{[\erf;1]}\left( -\frac{2t}{\Delta t}+1\right),
    \end{align}
where the time-dependent $g(t)$ utilises the compact mask function $\Gamma_{[\erf;1]}$ to smoothly accelerate from $g(0) = 0$ to $g(\Delta t) = 1$ at $\Delta t = 4$, with steady uniform flow boundary conditions ${g(t>\Delta t) = 1}$ until the end time $t_{end} = 10$.
The inner coordinate singularity at $r = 0$ is avoided by prescribing rotating no slip inner boundaries within the mask at $r = R_1 = 1/10$,
    \begin{align}
	u(R_1) &= 0, &    v(R_1) &= R_1 \Omega(t), & \text{where}&& \Omega(t) &= \sin\br{\frac{2t}{\pi}},
    \end{align}
where the cylinder angular velocity $\Omega$ is sinusoidal.
The reference simulation replaces the volume penalisation with the rotating no-slip boundary conditions applied at $r= R_1 = 1$.    
    
{\textbf{Force calculations}} --- 
The force $F_0$ and torque $T_0$ on a no slip boundary are defined by surface integrals of the tangential and normal components of the stress tensor
    \begin{align}
    F_{0,x} &= \int_{0}^{2\pi} \left[ r\left( \cos\theta (-P+\frac{2}{\Re}\partial_r u ) - \sin \theta \frac{1}{\Re r}(\partial_\theta u + r \partial_r v - v)\right)\right]_{r=R_1}\!\!d\theta,\\
    F_{0,y} &= \int_{0}^{2\pi} \left[ r\left( \sin\theta (-P+\frac{2}{\Re}\partial_r u ) + \cos \theta \frac{1}{\Re r}(\partial_\theta u + r \partial_r v - v)\right)\right]_{r=R_1}\!\!d\theta,\\
	T_{0} &= \int_{0}^{2\pi} \left[ \frac{r}{\Re }(\partial_\theta u + r \partial_r v - v)\right]_{r=R_1} \!\!d\theta.
    \end{align}
The penalised force $F$ and torque $T$ calculations of \cref{eq:standard-force-calculation} including correction terms are
    \begin{align}
    F_{x} &= \int_0^{2\pi}\!\! \int_{R_1}^{R_2} \frac{\Gamma}{\eta} \left( \cos\theta \,u - \sin \theta\, (v - r\Omega)\right) r \, dr \, d\theta + F_{0,x},\\
    F_{y} &= \int_0^{2\pi}\!\! \int_{R_1}^{R_2} \frac{\Gamma}{\eta} \left( \sin\theta \,u + \cos \theta\, (v - r\Omega)\right) r \, dr \, d\theta + F_{0,y},\\
	T &= \int_0^{2\pi}\!\! \int_{R_1}^{R_2} \frac{\Gamma}{\eta} (v - r\Omega) \, r^2 \,dr\,d\theta + \int_0^{2\pi}\!\!\int_{R_1}^1 \dot{\Omega} \, r^3 \, dr \, d\theta + T_0.
    \end{align}
Azimuthal symmetry means the acceleration correction appears only for the torque.

{\textbf{Volume-penalty parameters}}---
Three different types of penalty mask are chosen,
\begin{itemize}
	\item $\Gamma_{0,0}$: Standard discontinuous mask $\ell = 0, \delta = 0$,
	\item $\Gamma_{\ell,0}$: Shifted discontinuous mask $\ell = \eps, \delta = 0$,
	\item $\Gamma_{0,\delta}$: Zero-shift smoothed mask $\ell = 0, \delta = 3.80171928 \eps$,
\end{itemize}
which are all defined as various limits of a smooth compactified error function
    \begin{align}
    \Gamma_{\ell,\delta}(r) = \Ghat_{[\erf;1]}\br{\tfrac{r - 1 - \ell}{\delta}}.
    \end{align}
We then examine these mask families for $\eta_1 = 10^{-2}$, $\eta_2 = 10^{-3}$, and $\eta_3 = 10^{-4}$.

\subsection{Numerical method and parameters}

The initial value problem is simulated using the Dedalus framework \cite{BurnsDedalusFlexibleFramework2020} \footnote{Visit the Dedalus home page at \url{http://dedalus-project.org/}.}.
Dedalus is a general-purpose framework for solving partial differential equations using spectral methods.
Dedalus allows users to construct domains from the direct product of spectral series including Fourier Series, Chebyshev polynomial series, and continuous segments of Chebyshev polynomials (a {compound} Chebyshev basis), and to enter systems of PDEs in plain text.
The terms on the left-hand side (LHS) of the equations must be linear and are discretised into sparse and banded matrices acting on the spectral coefficients of the solution.
Boundary conditions are enforced at the endpoints of the Chebyshev bases with the use of a preconditioned Chebyshev-tau method.
The equations can be evolved using a range of mixed implicit-explicit timesteppers, with the LHS terms being integrated implicitly, and the right-hand side (RHS) terms integrated explicitly.
The framework is written in the Python programming language but utilises compiled libraries for performance and automatically parallelises the solvers using MPI, enabling for rapid prototyping and comparisons between models as well as efficient simulations at scale.

The test and reference solutions use well-resolved cylindrical grids in Dedalus to calculate the mathematical error of the penalised problem.
We use a Fourier basis of 512 modes for the azimuthal direction and Chebyshev bases in the radial direction.
A single Chebyshev basis with 256 modes is used for the reference solution, over the interval $I = [1, R_2]$. No-slip boundary conditions are enforced at $r=1$.
A compound Chebyshev basis is used for the penalised simulations, with additional segments to resolve the mask and the interior of the cylinder.
We partition over the intervals $I_1 = [R_1,1+\ell - \delta], I_2 = [1+\ell-\delta,1+\ell+\delta]$, and $I_3 = [1+\ell+\delta,R_2]$, where $R_1 = 0.1$ and $R_2 = 10$.
Each interval is discretised using 256 modes, 64 modes, and 256 modes respectively.
In the singular case $\delta = 0$ the middle interval is neglected.

The system is timestepped using a third order accurate backward difference scheme.
Other timesteppers (including Crank-Nicolson Adams-Bashforth and Runge-Kutta schemes) performed poorly with the time varying velocity boundary conditions for $t < 4$.
The scheme was validated as converged by observing spectral accuracy ($\O(10^{-10})$ error) for the pressure equation and boundary conditions at all times during the simulation.
The penalty term is timestepped explicitly (for accuracy reasons \cite{KolomenskiyFourierSpectralMethod2009}) with the stability constraint suggesting our choice of timestep $5\times 10^{-5}$.
This system is then timestepped to $t = 10$ for all simulations.
Minimal working examples of all simulations in this paper are available on GitHub at \cite{HesterMinimalworkingexamples2019}.

    \begin{table}[h!]
	\caption{Table of fixed cylinder simulation parameters.}
	\begin{center}
	\begin{tabular}{l|c|l}
	Parameter & Symbol & Value \Tstrut\Bstrut \\\hline 
	Reynolds number & $\Re$ & 200\\
	Inner domain radius & $R_1$ & 0.1\\
	Outer domain radius & $R_2$ & 10\\
	Final time & $t_{end}$ & 10\\
	Azimuthal modes & $N_\theta$ & 512
	\end{tabular}
	\end{center}
	\label{tab:fixed-cylinder-params}
    \end{table}

    \begin{table}[h!]
	\caption{Table of varying cylinder simulation parameters.}
	\begin{center}
	\begin{tabular}{c|c|c|c|c}
	Simulation & $\ell$ & $\delta$ & $\eta$ & $\eps$\Tstrut\Bstrut\\\hline
	$(0,0)$ & 0 & 0 & \num{1e-2}	& \num{7.07e-3}\\
	$(0,1)$ & 0 & 0 & \num{1e-3}	& \num{2.23e-3}\\	
	$(0,2)$ & 0 & 0 & \num{1e-4}	& \num{7.07e-4}\\	
	$(1,0)$ & \num{7.07e-2} & 0 & \num{1e-2}	& \num{7.07e-2}\\
	$(1,1)$ & \num{7.07e-3} & 0 & \num{1e-3}	& \num{2.23e-3}\\
	$(1,2)$ & \num{7.07e-4} & 0 & \num{1e-4}	& \num{7.07e-4}\\
	$(2,0)$ & 0 & \num{2.69e-2} & \num{1e-2}	& \num{7.07e-2}\\
	$(2,1)$ & 0 & \num{8.50e-3} & \num{1e-3}	& \num{2.23e-3}\\
	$(2,2)$ & 0 & \num{2.69e-3} & \num{1e-4}	& \num{7.07e-4}\\
	\end{tabular}
	\end{center}
	\label{tab:varying-cylinder-params}
    \end{table}

\subsection{Richardson extrapolation of model error}
{Richardson extrapolation} is a general technique to accelerate convergence of sequences with known asymptotic behaviour.
For a quantity $X_i,X_j$ derived from simulations with penalty parameters $\eta_i$ and $\eta_j$ we can calculate an {extrapolated} quantity $\overline{X}_{ij}$ as
    \begin{align}
    \overline{X}_{ij} \equiv \frac{X_i/\eta_i - X_j/\eta_j}{1/\eta_i - 1/\eta_j}.
    \end{align}
If $X_i,X_j$ obey leading order behaviour proportional to $\eta$, the extrapolation $\overline{X}_{ij}$ will cancel this error.
One can thus achieve $\O(\eta^2)$ accuracy in practice, for the small added cost of a lower resolution simulation at larger $\eta$.

\subsection{Summary of errors}
\begin{figure}[ht!]
	\centering
  \includegraphics[width=.8\linewidth]{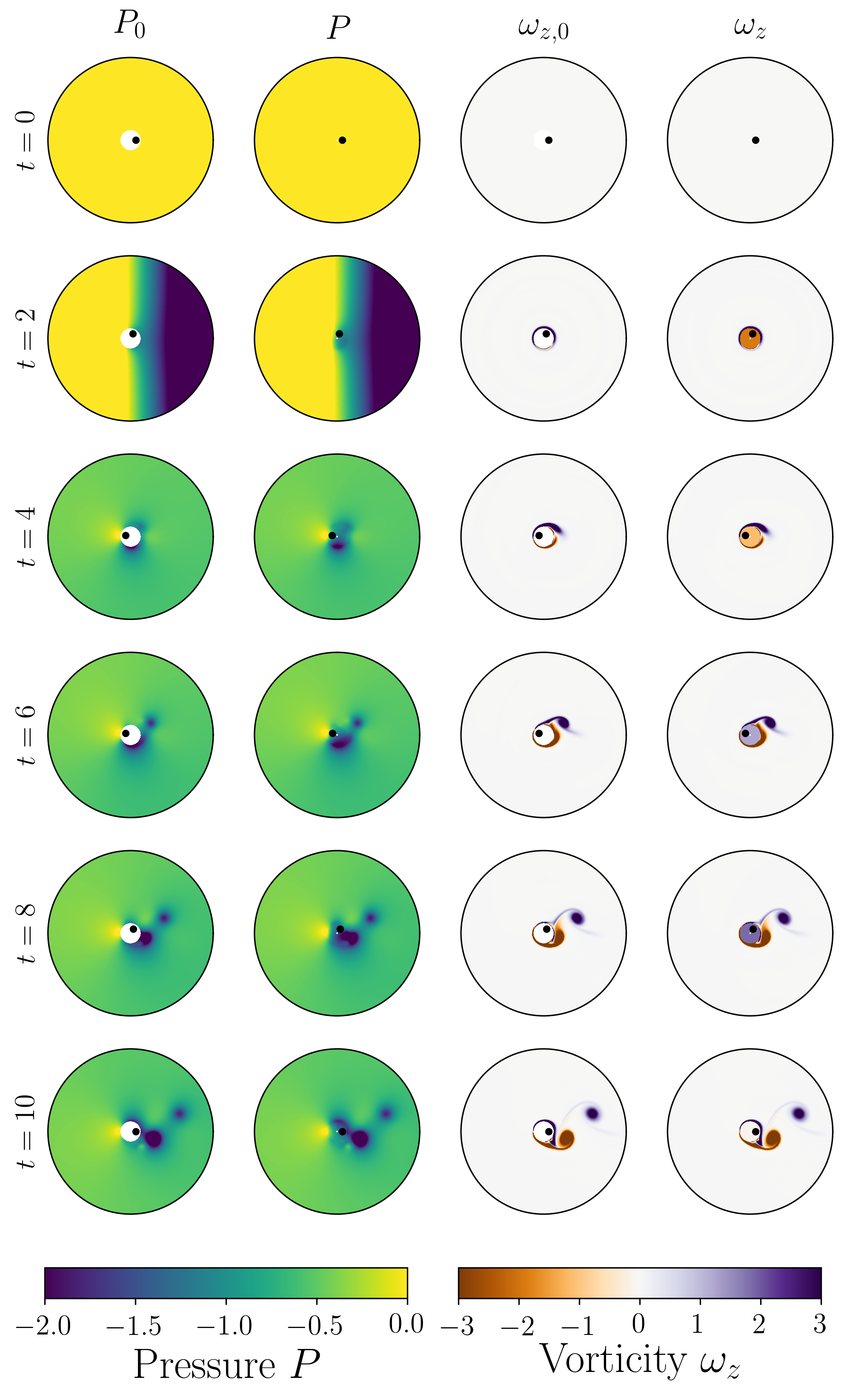}
  \caption{Time snapshots of pressure $P=p-(u^2+v^2)/2$, and vorticity $\omega_z=q/r$ for the reference and (2,2) volume-penalty cylinder simulations. The interior black dot shows the orientation of the rotating cylinder. Rotation of the volume penalised cylinder results in non-zero vorticity within the object (column 4). The pressure and vorticity fields are indistinguishable outside the object however.}
  \label{fig:oscillation-comparison-plots}
\end{figure}

\Cref{fig:oscillation-comparison-plots} compares the pressure $P$ and vorticity $\omega_z$ fields for the reference simulation and the volume-penalty simulation $(2,2)$.
Snapshots are shown every two seconds until time $t = 10$, and a black dot illustrates the rotation of the cylinder.
The pressure color scale is calibrated to the steady state uniform flow after $t = 4$,
and the figures are plotted to $r = 8$ to capture the vortex evolution. 
After $t=4$ the outer boundaries are steady and higher pressures are seen at the streamwise side of the cylinder, with low pressures around vortices and the sides of the cylinder.
The pressure varies {within} the volume penalised cylinder, requiring the inner stress correction term in the force and torque calculations.
Vortex shedding is observed immediately due to the rotation of the cylinder and the inherent instability of the von Karman vortex street.

\begin{figure}[ht!]
\centering
  \includegraphics[width=.7\linewidth]{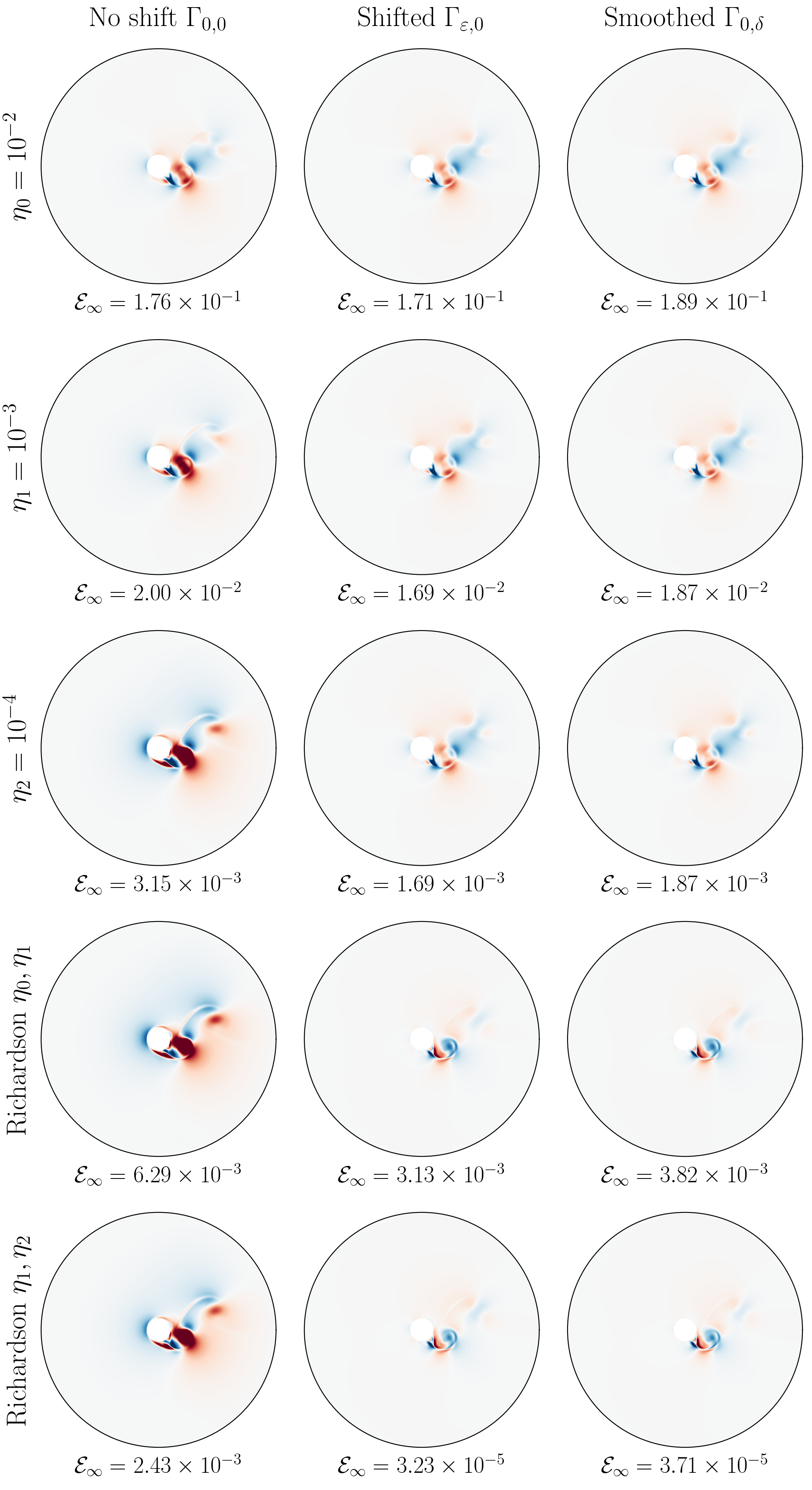}
  \caption{Error pattern between penalised and reference radial velocities $u - u_0$. Plots are made at $t = 10$ for the three different mask choices, at three choices of $\eta$. Each color map is normalised to the maximum error of the velocity field $\E_\infty$. The spatial pattern of the error is almost identical for all optimised masks, and differs only in the magnitude as $\eta$ is decreased. Two Richardson extrapolated fields are calculated for each mask type, which significantly improve the error for the optimised masks, but fails for the unoptimised mask.}
    \label{fig:oscillation-u-errors}
\end{figure}

To quantify the error \cref{fig:oscillation-u-errors} plots the pointwise error $u - u_0$ in the radial velocity field at $t=10$ for the different masks ($\Gamma_{0,0}, \Gamma_{\ell,0}, \Gamma_{0,\delta}$) and penalty parameters $\eta = \eta_0, \eta_1, \eta_2$.
The colormap range is normalised to half the maximum local error $\E_\infty$, printed below each figure.
Optimal adjustments give identical {spatial profiles} of the error (the two right columns), with the amplitude decreasing as $\O(\eta)$.
Standard discontinuous masks lack this property, with varying spatial pattern and reduced convergence in $\eta$.
This is strong empirical evidence for the validity of the asymptotic expansions of \cref{sec:asymptotics}.

This motivates the extrapolated calculations in \cref{fig:oscillation-u-errors}.
Extrapolations of consecutive values of $\eta$ were performed for each mask.
No improvement is observed for the standard discontinuous mask, as expected from the $\O(\eps)$ displacement length error.
However extrapolation is very effective for {optimised} masks.
Extrapolation using $\eta_0 = 10^{-2}$ and $\eta_1^{-3}$ performs similarly to $\eta_2 = 10^{-4}$, and extrapolation from $\eta_1$ and $\eta_2$ is almost two orders of magnitude more accurate than $\eta_2$ alone.
This shows $\O(\eta^2)$ convergence for Richardson extrapolations of optimal masks, a massive performance boost.
By running a smaller second simulation in tandem with the first, it is possible to improve convergence to $\O(\eta^2) \propto \O(\eps^4)$, implying only twice the effort to halve the error in 2D, and slightly more in 3D, {vastly} more efficient than the 16 and 32-fold increases for uncorrected masks.
This extrapolation is performed post simulation, but could be calculated {during} a run.
The Lyapunov time of chaotic fluid flows would require reinitializing each simulation over this time scale to prevent desynchronisation.

\begin{figure}[ht!]
	\centering
  \includegraphics[width=.9\linewidth]{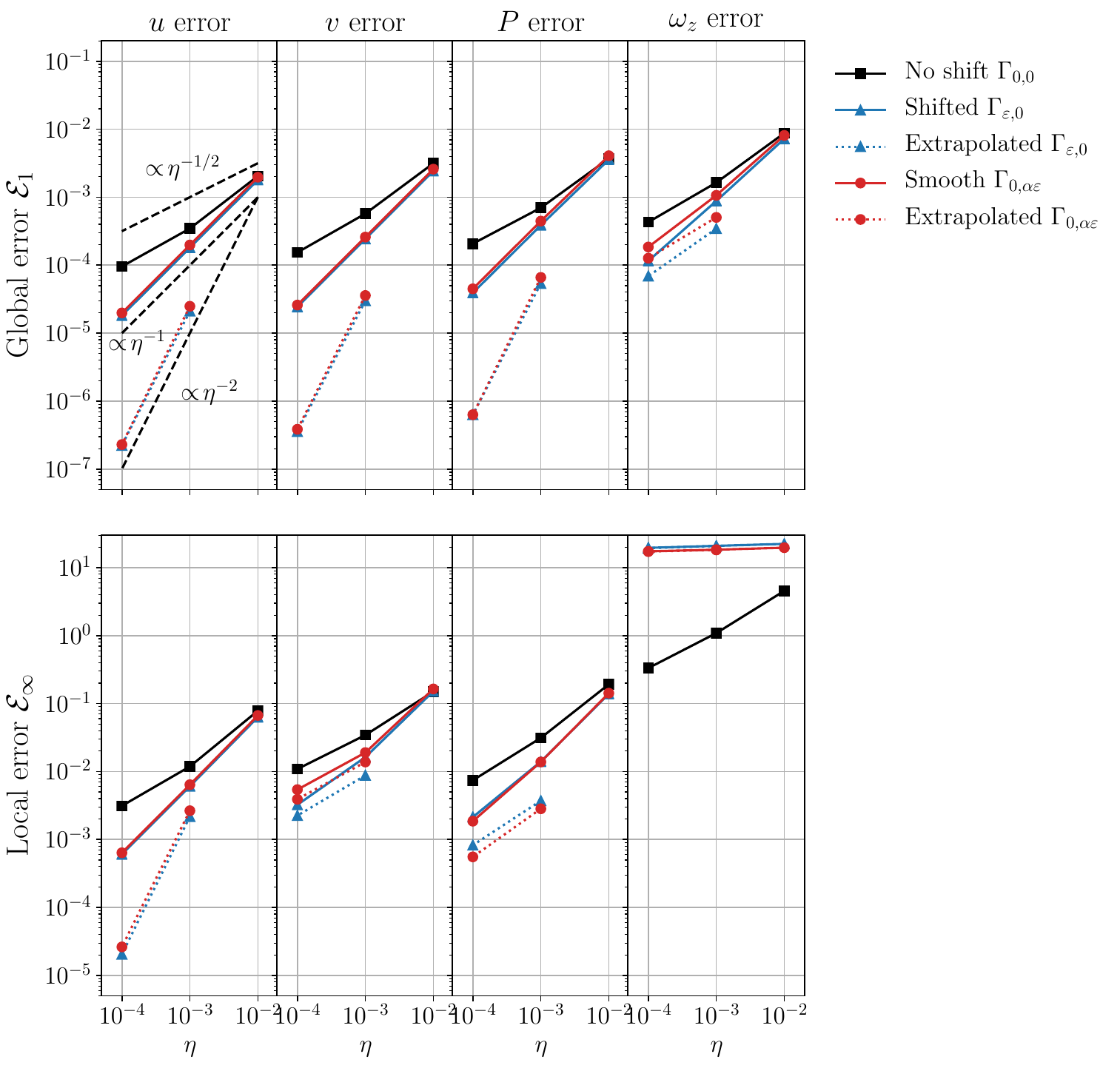}
  \caption{Plots of the global $\E_1$ and local error $\E_\infty$ for $u,v,P,$ and $\omega_z$ as a function of $\eta$, for three mask choices ($\Gamma_{0,0}$ (black), $\Gamma_{\eps,0}$ (blue), $\Gamma_{0,\delta}$ (red)). Errors were calculated for $1<r<10$ at $t = 10$. Extrapolated calculations are dotted. Different convergence behaviours are apparent for different variables and different masks.}
    \label{fig:oscillation-variable-errors}
\end{figure}

We then calculate the global $\E_1$ and local $\E_\infty$ error norms at the final timestep for the radial velocity $u$, the azimuthal velocity $v$, the true pressure $P$, and the true vorticity $\omega_z$ in \cref{fig:oscillation-variable-errors}.
We see the standard discontinuous mask $\Gamma_{0,0}$ (black) performs significantly worse than optimised masks, giving the predicted $\O(\eta^{1/2})$ convergence rate.
Optimised masks instead show $\O(\eta)$ convergence, and extrapolated predictions achieve $\O(\eta^2)$ error (based off the smallest $\eta$ used).
However, for the weakest penalisation $\eta_0 = \num{1e-2}$, all approaches are equally accurate;
In the intermediate damping regime, the dominant error is not due to the displacement length (which optimised masks lack), but {time scale} errors due to  pressure gradients and acceleration of the penalty mask, as seen for stagnation point flow in \cref{sec:inner}.
Similar behaviour for the global error is observed for the azimuthal velocity $v$ and pressure $P$.

The differing behaviour of the global and local errors for each variable stems from the differentiability of the true solutions at the boundary, and the localisation of error achieved by optimal masks.
The radial velocity is $C^1$ as it and its derivative are zero at the boundary. 
This is uniformly approximated by penalised solutions, giving commensurate global $\E_1$ and local $\E_\infty$ errors.
However the tangential velocity $v_0$ and pressure $p_0$ are $C^0$ with an $\O(1)$ jump in the derivative at the boundary (e.g. \cref{fig:inner-diagrams}).
This is less regular than the volume penalised solution, causing $\O(\eps)$ disagreement in the boundary region.
Calibrating the mask localises this error to the boundary, implying maximum error of $\O(\eps)$, but global error of only $\O(\eps^2)$.
The derivative (vorticity) is {discontinuous} at the boundary, causing $\O(1)$ disagreement in the interior of the mask near the boundary.
Optimal masks cause $\O(1)$ disagreement in $\omega_z$ in an $\O(\eps)$ boundary region and $\O(\eps^2)$ elsewhere, while standard masks cause $\O(\eps)$ error everywhere.

\begin{figure}[ht!]
	\centering
  \includegraphics[width=.85\linewidth]{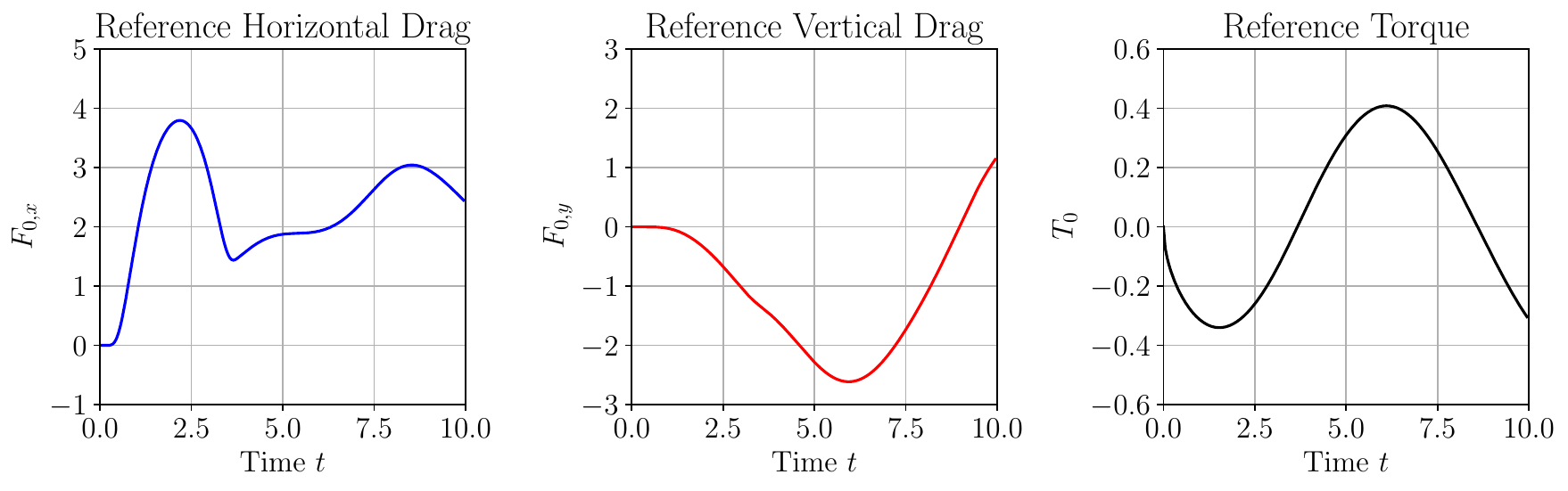}
  \caption{Time series of the total horizontal drag $F_{0,x}$, total vertical drag $F_{0,y}$, and total torque $T_0$ for the reference cylinder simulation at $\Re = 200$.}
  \label{fig:oscillation-reference-drags}
\end{figure}

We plot time series of the reference horizontal drag $F_{0,x}$, vertical lift $F_{0,y}$, and reference torque $T_0$ in \cref{fig:oscillation-reference-drags}.
Heightened drag is observed during the peak fluid acceleration ($t < \Delta t = 4$), corresponding to the large induced pressure gradient observed at $t = 2$ in \cref{fig:oscillation-reference-drags}.
A significant lift force is also generated during the simulation, as expected from the Magnus effect \cite{MagnusUeberAhweichungGeschosse1853}.

\begin{figure}[ht!]
	\centering
  \includegraphics[width=.95\linewidth]{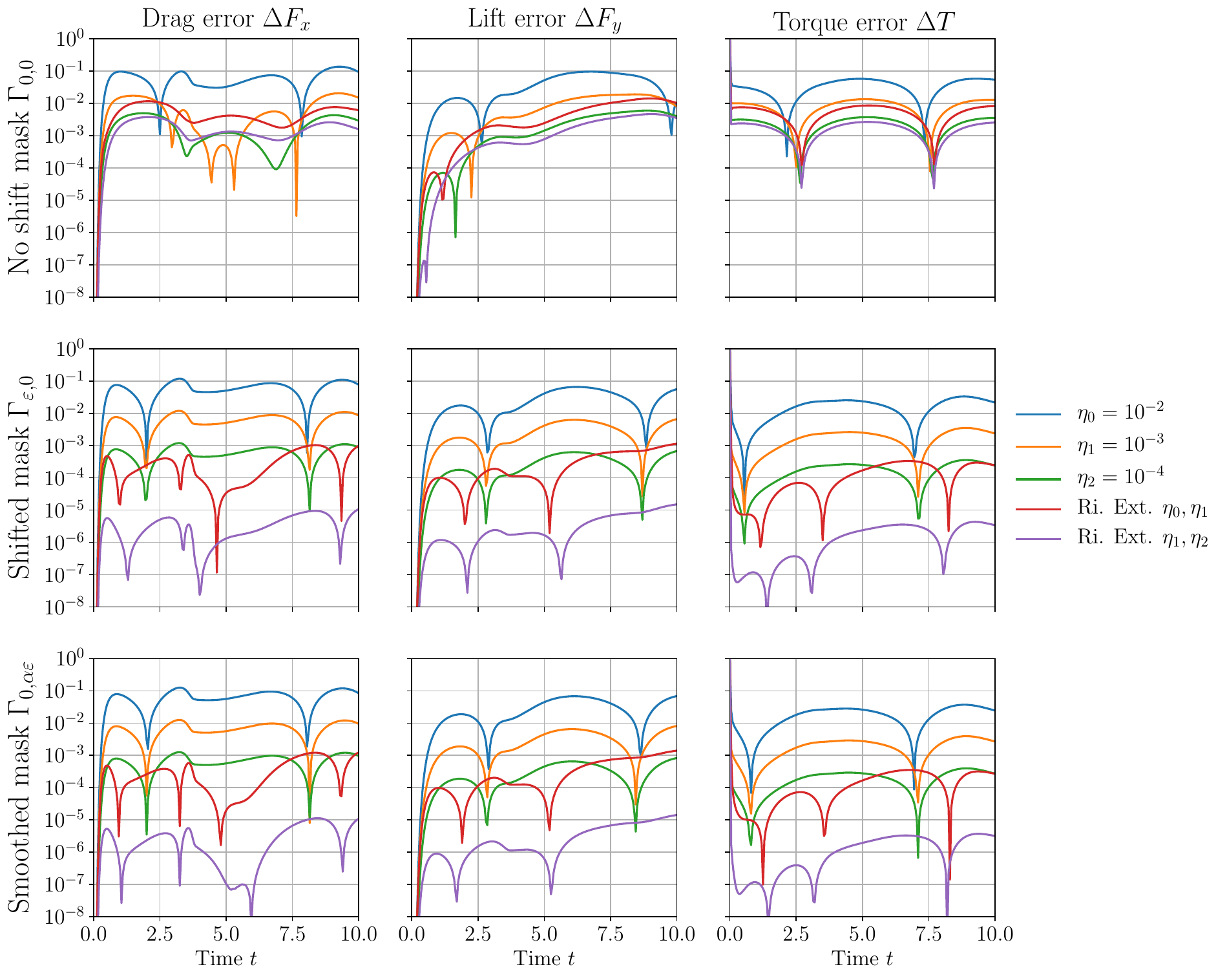}
  \caption{Time series of the drag error $\Delta F_x$, lift error $\Delta F_y$, and torque error $\Delta T$ at $\Re = 200$ for standard discontinuous $\Gamma_{0,0}$, optimal shifted $\Gamma_{\eps,0}$, and optimal smooth $\Gamma_{0,\alpha\eps}$ masks for different choices of damping $\eta = 10^{-2},10^{-3},10^{-4}$, and extrapolated errors.}
  \label{fig:oscillation-physical-errors}
\end{figure}

\Cref{fig:oscillation-physical-errors} compares the time series of the horizontal drag error $\Delta F_x$, vertical drag error $\Delta F_y$, and torque error $\Delta T$ for each choice of mask $\Gamma_{0,0},\Gamma_{\ell,0}$ and $\Gamma_{0,\delta}$ and choice of damping $\eta=\eta_0,\eta_1,\eta_2$.
All the masks improve similarly in the transition from $\eta_1 = 10^{-2}$ to $\eta_2 = 10^{-3}$, in keeping with the intermediate damping regime.
However, beyond $\eta = 10^{-2}$ the unshifted discontinuous mask $\Gamma_{0,0}$ transitions from $\O(\Re \eps^2)$ to $\O(\eps)$ convergence,
whereas optimised masks maintain consistent $\O(\eta)$ convergence.
The temporal profile of the error is largely identical for the corrected masks, but inconsistent for the unshifted discontinuous mask.
Applying the extrapolation procedure to the drag, lift and torque calculations, and we find similar performance improvements.
No improvement occurs for the standard mask, while the adjusted masks each show $\O(\eta^2)$ convergence.
This is an extraordinary improvement in the accuracy of the volume-penalty method, attaining drag errors of order $10^{-5}$ using extrapolation at $\eta = 10^{-4}$.

\section{Conclusions}
\label{sec:conclusion}

In this paper we examine and improve the volume-penalty method, a model of objects in fluids which approximates standard no-slip boundaries.
A key observation is that damping induces both a penalisation {time} scale $\eta$, and {length} scale $\eps$, related by the Reynolds number $\eta = \Re \, \eps^2$.
By using multiple-scales matched-asymptotics to second order in $\eps$, and convenient signed distance boundary layer coordinates, we show that the leading order model error in the fluid is in general $\O(\eps)$ due to the {displacement length} of general mask functions.

We then show this displacement error can be eliminated by correctly {shifting} and {smoothing} the mask near the boundary, promoting the accuracy to $\O(\eta)$. 
This calculation can be done once efficiently for a particular mask function using a {Riccati} transformation.
This improvement in model error translates to an improvement in total error from first to second order numerical convergence.

The second order asymptotics also reveals two damping regimes for large $\Re$. 
{Intermediate} damping $1 > \eta > \eps > \Re^{-1}$ implies the `second order' time scale error dominates the `first order' displacement length error,
while {strong} damping $\Re^{-1} > \eps > \eta$ shows the displacement length error dominates and mask corrections must be used to achieve $\O(\eta)$ accuracy.
We also derive simple force and torque calculations which are proportional to the accuracy of the flow. 

We validate each of these findings in Poiseuille flow, viscous stagnation point flow, and accelerating flow past a rotating cylinder at $\Re = 200$.
We finally show that in the strong damping regime, Richardson extrapolation can also be used with optimised masks to boost model error convergence to $\O(\eta^2)$.

\section*{Acknowledgments}
Eric Hester acknowledges support from The University of Sydney through the Postgraduate Teaching Fellowship, and the Phillip Hofflin International Research Travel Scholarship.
Geoffrey Vasil acknowledges support from the Australian Research Council, project number DE140101960.
The authors would also like to acknowledge the thoughtful input from many discussions with Daniel Lecoanet, Christopher Lustri, Benjamin Favier, Laurent Duchamin, and Kai Schneider.

\bibliographystyle{abbrv}
\bibliography{vp-references}

\begin{thebibliography}{10}

\bibitem{AndersonDiffuseInterfaceMethodsFluid1998}
D.~M. Anderson, G.~B. McFadden, and A.~A. Wheeler.
\newblock Diffuse-{{Interface Methods}} in {{Fluid Mechanics}}.
\newblock {\em Annual Review of Fluid Mechanics}, 30(1):139--165, 1998.

\bibitem{AngotAnalysisSingularPerturbations1999}
P.~Angot.
\newblock Analysis of singular perturbations on the {{Brinkman}} problem for
  fictitious domain models of viscous flows.
\newblock {\em Mathematical Methods in the Applied Sciences},
  22(16):1395--1412, Nov. 1999.

\bibitem{AngotPenalizationMethodTake1999}
P.~Angot, C.-H. Bruneau, and P.~Fabrie.
\newblock A penalization method to take into account obstacles in
  incompressible viscous flows.
\newblock {\em Numerische Mathematik}, 81(4):497--520, Feb. 1999.

\bibitem{BlackburnSemtexSpectralElement2019}
H.~M. Blackburn, D.~Lee, T.~Albrecht, and J.~Singh.
\newblock Semtex: {{A}} spectral element\textendash{{Fourier}} solver for the
  incompressible {{Navier}}\textendash{{Stokes}} equations in cylindrical or
  {{Cartesian}} coordinates.
\newblock {\em Computer Physics Communications}, 245:106804, Dec. 2019.

\bibitem{BoettingerPhaseFieldSimulationSolidification2002}
W.~J. Boettinger, J.~A. Warren, C.~Beckermann, and A.~Karma.
\newblock Phase-{{Field Simulation}} of {{Solidification}}.
\newblock {\em Annual Review of Materials Research}, 32(1):163--194, 2002.

\bibitem{BoydFourierEmbeddedDomain2005}
J.~P. Boyd.
\newblock Fourier embedded domain methods: Extending a function defined on an
  irregular region to a rectangle so that the extension is spatially periodic
  and {{C}}{$\infty$}.
\newblock {\em Applied Mathematics and Computation}, 161(2):591--597, Feb.
  2005.

\bibitem{BriscoliniDevelopmentMaskMethod1989}
M.~Briscolini and P.~Santangelo.
\newblock Development of the mask method for incompressible unsteady flows.
\newblock {\em Journal of Computational Physics}, 84(1):57--75, Sept. 1989.

\bibitem{BurgerAnalysisDiffuseDomain2017}
M.~Burger, O.~L. Elvetun, and M.~Schlottbom.
\newblock Analysis of the {{Diffuse Domain Method}} for {{Second Order Elliptic
  Boundary Value Problems}}.
\newblock {\em Foundations of Computational Mathematics}, 17(3):627--674, June
  2017.

\bibitem{BurnsDedalusFlexibleFramework2020}
K.~J. Burns, G.~M. Vasil, J.~S. Oishi, D.~Lecoanet, and B.~P. Brown.
\newblock Dedalus: {{A}} flexible framework for numerical simulations with
  spectral methods.
\newblock {\em Physical Review Research}, 2(2):023068, Apr. 2020.

\bibitem{CaginalpAnalysisPhaseField1986}
G.~Caginalp.
\newblock An analysis of a phase field model of a free boundary.
\newblock {\em Archive for Rational Mechanics and Analysis}, 92(3):205--245,
  Sept. 1986.

\bibitem{CaginalpStefanHeleShawType1989}
G.~Caginalp.
\newblock Stefan and {{Hele}}-{{Shaw}} type models as asymptotic limits of the
  phase-field equations.
\newblock {\em Physical Review A}, 39(11):5887--5896, June 1989.

\bibitem{CahnFreeEnergyNonuniform1958}
J.~W. Cahn and J.~E. Hilliard.
\newblock Free {{Energy}} of a {{Nonuniform System}}. {{I}}. {{Interfacial Free
  Energy}}.
\newblock {\em The Journal of Chemical Physics}, 28(2):258--267, Feb. 1958.

\bibitem{CantwellNektarOpensourceSpectral2015}
C.~D. Cantwell, D.~Moxey, A.~Comerford, A.~Bolis, G.~Rocco, G.~Mengaldo,
  D.~De~Grazia, S.~Yakovlev, J.~E. Lombard, D.~Ekelschot, B.~Jordi, H.~Xu,
  Y.~Mohamied, C.~Eskilsson, B.~Nelson, P.~Vos, C.~Biotto, R.~M. Kirby, and
  S.~J. Sherwin.
\newblock Nektar++: {{An}} open-source spectral/hp element framework.
\newblock {\em Computer Physics Communications}, 192:205--219, July 2015.

\bibitem{CarbouBoundaryLayerPenalization2003}
G.~Carbou and P.~Fabrie.
\newblock Boundary layer for a penalization method for viscous incompressible
  flow.
\newblock {\em Advances in Differential Equations}, 8(12):1453--1480, 2003.

\bibitem{ChenRapidlyConvergingPhase2006}
X.~Chen, G.~Caginalp, and C.~Eck.
\newblock A rapidly converging phase field model.
\newblock {\em Discrete and Continuous Dynamical Systems}, 15(4):1017--1034,
  May 2006.

\bibitem{ColeQuasiLinearParabolicEquation1951}
J.~D. Cole.
\newblock On a {{Quasi}}-{{Linear Parabolic Equation Occurring In
  Aerodynamics}}.
\newblock {\em Quarterly of Applied Mathematics}, 9(3):225--236, 1951.

\bibitem{DevilleHighorderMethodsIncompressible2002}
M.~Deville, P.~Fischer, P.~Fischer, and E.~Mund.
\newblock {\em High-Order Methods for Incompressible Fluid Flow}.
\newblock Cambridge Monographs on Applied and Computational Mathematics.
  {Cambridge University Press}, 2002.

\bibitem{EngelsNumericalSimulationFluid2015}
T.~Engels, D.~Kolomenskiy, K.~Schneider, and J.~Sesterhenn.
\newblock Numerical simulation of fluid\textendash structure interaction with
  the volume penalization method.
\newblock {\em Journal of Computational Physics}, 281(Supplement C):96--115,
  Jan. 2015.

\bibitem{EngelsFluSINovelParallel2016}
T.~Engels, D.~Kolomenskiy, K.~Schneider, and J.~Sesterhenn.
\newblock {{FluSI}}: {{A Novel Parallel Simulation Tool}} for {{Flapping Insect
  Flight Using}} a {{Fourier Method}} with {{Volume Penalization}}.
\newblock {\em SIAM Journal on Scientific Computing}, 38(5):S3--S24, Jan. 2016.

\bibitem{FischerNek5000WebPage2008}
P.~F. Fischer, S.~G. Kerkemeier, and J.~W. Lottes.
\newblock Nek5000 {{Web}} page.
\newblock http://nek5000.mcs.anl.gov, 2008.

\bibitem{FranzNoteConvergenceAnalysis2012}
S.~Franz, H.-G. Roos, R.~G{\"a}rtner, and A.~Voigt.
\newblock A {{Note}} on the {{Convergence Analysis}} of a {{Diffuse}}-domain
  {{Approach}}.
\newblock {\em Computational Methods in Applied Mathematics}, 12(2), 2012.

\bibitem{GibbsEquilibriumHeterogeneousSubstances1878}
J.~W. Gibbs.
\newblock On the equilibrium of heterogeneous substances.
\newblock {\em American Journal of Science}, Series 3 Vol. 16(96):441--458,
  Dec. 1878.

\bibitem{GlowinskiFictitiousDomainApproach2001}
R.~Glowinski, T.~W. Pan, T.~I. Hesla, D.~D. Joseph, and J.~P{\'e}riaux.
\newblock A {{Fictitious Domain Approach}} to the {{Direct Numerical
  Simulation}} of {{Incompressible Viscous Flow}} past {{Moving Rigid Bodies}}:
  {{Application}} to {{Particulate Flow}}.
\newblock {\em Journal of Computational Physics}, 169(2):363--426, May 2001.

\bibitem{GlowinskiFictitiousDomainMethod1994}
R.~Glowinski, T.-W. Pan, and J.~Periaux.
\newblock A fictitious domain method for {{Dirichlet}} problem and
  applications.
\newblock {\em Computer Methods in Applied Mechanics and Engineering},
  111(3):283--303, Jan. 1994.

\bibitem{GoldsteinModelingNoSlipFlow1993}
D.~Goldstein, R.~Handler, and L.~Sirovich.
\newblock Modeling a {{No}}-{{Slip Flow Boundary}} with an {{External Force
  Field}}.
\newblock {\em Journal of Computational Physics}, 105(2):354--366, Apr. 1993.

\bibitem{HesterMinimalworkingexamples2019}
E.~W. Hester.
\newblock Minimal working examples for volume penalized boundary conditions in
  {{Dedalus}}., May 2019.

\bibitem{HohenbergTheoryDynamicCritical1977}
P.~C. Hohenberg and B.~I. Halperin.
\newblock Theory of dynamic critical phenomena.
\newblock {\em Reviews of Modern Physics}, 49(3):435--479, July 1977.

\bibitem{HopfPartialDifferentialEquation1950}
E.~Hopf.
\newblock The partial differential equation ut + uux = {$M$}xx.
\newblock {\em Communications on Pure and Applied Mathematics}, 3(3):201--230,
  1950.

\bibitem{IaccarinoImmersedBoundaryTechnique2003}
G.~Iaccarino and R.~Verzicco.
\newblock Immersed boundary technique for turbulent flow simulations.
\newblock {\em Applied Mechanics Reviews}, 56(3):331--347, May 2003.

\bibitem{JacqminCalculationTwoPhaseNavier1999}
D.~Jacqmin.
\newblock Calculation of {{Two}}-{{Phase Navier}}\textendash{{Stokes Flows
  Using Phase}}-{{Field Modeling}}.
\newblock {\em Journal of Computational Physics}, 155(1):96--127, Oct. 1999.

\bibitem{KadochVolumePenalizationMethod2012}
B.~Kadoch, D.~Kolomenskiy, P.~Angot, and K.~Schneider.
\newblock A volume penalization method for incompressible flows and scalar
  advection\textendash diffusion with moving obstacles.
\newblock {\em Journal of Computational Physics}, 231(12):4365--4383, June
  2012.

\bibitem{KarmaPhasefieldMethodComputationally1996}
A.~Karma and W.-J. Rappel.
\newblock Phase-field method for computationally efficient modeling of
  solidification with arbitrary interface kinetics.
\newblock {\em Physical Review E}, 53(4):R3017--R3020, Apr. 1996.

\bibitem{KarniadakisSpectralHpElement2005}
G.~Karniadakis and S.~Sherwin.
\newblock {\em Spectral/Hp Element Methods for Computational Fluid Dynamics}.
\newblock {Oxford University Press}, 2nd edition edition, 2005.

\bibitem{KevlahanComputationTurbulentFlow2001}
N.~K.~R. Kevlahan and J.-M. Ghidaglia.
\newblock Computation of turbulent flow past an array of cylinders using a
  spectral method with {{Brinkman}} penalization.
\newblock {\em European Journal of Mechanics - B/Fluids}, 20(3):333--350, May
  2001.

\bibitem{KhadraFictitiousDomainApproach2000}
K.~Khadra, P.~Angot, S.~Parneix, and J.-P. Caltagirone.
\newblock Fictitious domain approach for numerical modelling of
  {{Navier}}\textendash{{Stokes}} equations.
\newblock {\em International Journal for Numerical Methods in Fluids},
  34(8):651--684, Dec. 2000.

\bibitem{KolomenskiyAnalysisDiscretizationVolume2015}
D.~Kolomenskiy, R.~{Nguyen van yen}, and K.~Schneider.
\newblock Analysis and discretization of the volume penalized {{Laplace}}
  operator with {{Neumann}} boundary conditions.
\newblock {\em Applied Numerical Mathematics}, 95:238--249, Sept. 2015.

\bibitem{KolomenskiyFourierSpectralMethod2009}
D.~Kolomenskiy and K.~Schneider.
\newblock A {{Fourier}} spectral method for the {{Navier}}\textendash{{Stokes}}
  equations with volume penalization for moving solid obstacles.
\newblock {\em Journal of Computational Physics}, 228(16):5687--5709, Sept.
  2009.

\bibitem{LangerModelsPatternFormation1986}
J.~S. Langer.
\newblock Models of {{Pattern Formation}} in {{First}}-{{Order Phase
  Transitions}}.
\newblock In {\em Directions in {{Condensed Matter Physics}}}, pages 165--186.
  1986.

\bibitem{LervagAnalysisDiffusedomainMethod2015}
K.~Y. Lerv{\aa}g and J.~Lowengrub.
\newblock Analysis of the diffuse-domain method for solving {{PDEs}} in complex
  geometries.
\newblock {\em Communications in Mathematical Sciences}, 13(6):1473--1500,
  2015.

\bibitem{LiSolvingPDEsComplex2009}
X.~Li, J.~Lowengrub, A.~Ratz, and A.~Voigt.
\newblock Solving {{PDEs}} in complex geometries: A diffuse domain approach.
\newblock {\em Communications in mathematical sciences}, 7(1):81--107, Mar.
  2009.

\bibitem{MagnusUeberAhweichungGeschosse1853}
G.~Magnus.
\newblock Ueber die {{Ahweichung}} der {{Geschosse}}, und: {{Ueber}} eine
  auffallende {{Erscheinung}} bei rotirenden {{K\"orpern}}.
\newblock {\em Annalen der Physik}, 164(1):1--29, 1853.

\bibitem{MittalImmersedBoundaryMethods2005}
R.~Mittal and G.~Iaccarino.
\newblock Immersed {{Boundary Methods}}.
\newblock {\em Annual Review of Fluid Mechanics}, 37(1):239--261, 2005.

\bibitem{PeskinNumericalAnalysisBlood1977}
C.~S. Peskin.
\newblock Numerical analysis of blood flow in the heart.
\newblock {\em Journal of Computational Physics}, 25(3):220--252, Nov. 1977.

\bibitem{PeskinImmersedBoundaryMethod2002}
C.~S. Peskin.
\newblock The immersed boundary method.
\newblock {\em Acta Numerica}, 11:479--517, Jan. 2002.

\bibitem{PoulsenSmoothedBoundaryMethod2017}
S.~O. Poulsen and P.~W. Voorhees.
\newblock Smoothed {{Boundary Method}} for {{Diffusion}}-{{Related Partial
  Differential Equations}} in {{Complex Geometries}}.
\newblock {\em International Journal of Computational Methods}, 15(03):1850014,
  July 2017.

\bibitem{SakuraiVolumePenalizationInhomogeneous2019}
T.~Sakurai, K.~Yoshimatsu, N.~Okamoto, and K.~Schneider.
\newblock Volume penalization for inhomogeneous {{Neumann}} boundary conditions
  modeling scalar flux in complicated geometry.
\newblock {\em Journal of Computational Physics}, 390:452--469, Aug. 2019.

\bibitem{SchneiderImmersedBoundaryMethods2015}
K.~Schneider.
\newblock Immersed boundary methods for numerical simulation of confined fluid
  and plasma turbulence in complex geometries: A review.
\newblock {\em Journal of Plasma Physics}, 81(6), Dec. 2015.

\bibitem{SchneiderPseudospectralMethodVolume2011}
K.~Schneider, S.~Neffaa, and W.~J.~T. Bos.
\newblock A pseudo-spectral method with volume penalisation for
  magnetohydrodynamic turbulence in confined domains.
\newblock {\em Computer Physics Communications}, 182(1):2--7, Jan. 2011.

\bibitem{StefanUeberTheorieEisbildung1891}
J.~Stefan.
\newblock Ueber die {{Theorie}} der {{Eisbildung}}, insbesondere \"uber die
  {{Eisbildung}} im {{Polarmeere}}.
\newblock {\em Annalen der Physik}, 278(2):269--286, 1891.

\bibitem{UhlmannImmersedBoundaryMethod2005}
M.~Uhlmann.
\newblock An immersed boundary method with direct forcing for the simulation of
  particulate flows.
\newblock {\em Journal of Computational Physics}, 209(2):448--476, Nov. 2005.

\bibitem{VanderWaalsThermodynamicTheoryCapillarity1979}
J.~D. {van der Waals}.
\newblock The thermodynamic theory of capillarity under the hypothesis of a
  continuous variation of density.
\newblock {\em Journal of Statistical Physics}, 20(2):200--244, Feb. 1979.

\bibitem{YenApproximationLaplaceStokes2014}
R.~N. van Yen, D.~Kolomenskiy, and K.~Schneider.
\newblock Approximation of the {{Laplace}} and {{Stokes}} operators with
  {{Dirichlet}} boundary conditions through volume penalization: A spectral
  viewpoint.
\newblock {\em Numerische Mathematik}, 128(2):301--338, Oct. 2014.

\bibitem{YuHigherorderAccurateDiffusedomain2020}
F.~Yu, Z.~Guo, and J.~Lowengrub.
\newblock Higher-order accurate diffuse-domain methods for partial differential
  equations with {{Dirichlet}} boundary conditions in complex, evolving
  geometries.
\newblock {\em Journal of Computational Physics}, 406:109174, Apr. 2020.

\bibitem{YuDLMFDMethod2005}
Z.~Yu.
\newblock A {{DLM}}/{{FD}} method for fluid/flexible-body interactions.
\newblock {\em Journal of Computational Physics}, 207(1):1--27, July 2005.

\bibitem{YuDirectforcingFictitiousDomain2007}
Z.~Yu and X.~Shao.
\newblock A direct-forcing fictitious domain method for particulate flows.
\newblock {\em Journal of Computational Physics}, 227(1):292--314, Nov. 2007.

\end{thebibliography}

\end{document}